\documentclass[amstex, amysymb, 12pt]{article}
\usepackage{amsmath, amsfonts, amssymb, amsthm}
\usepackage{graphicx}

\newtheorem{Pa}{Paper}[section]
\newtheorem{Tm}[Pa]{{\bf Theorem}}
\newtheorem{La}[Pa]{{\bf Lemma}}

\newtheorem{Rk}[Pa]{{\bf Remark}}

\newtheorem{Dn}[Pa]{{\bf Definition}}

\newcommand{\nline}{{\vspace{\baselineskip}}}

\newcounter{caser}[section]
\newcommand{\Cs}
 {\addtocounter{caser}{1}\vspace{5mm}\par\noindent \underline{{\bf Case}  \arabic{caser}}:\space}
\begin{document}

\title{The Graphs for which the Maximum Multiplicity of an Eigenvalue
is Two}

\author{
Charles R. Johnson\footnote{Corresponding author;
Department of Mathematics, College of William and Mary, Williamsburg, VA  23187-8795, USA;
email: crjohnso@math.wm.edu; phone (office): 1-757-221-2014; fax (identify recipient): 1-757-221-7400}, 
\and Raphael Loewy\footnote{Department of Mathematics, Technion -- Israel Institute of Technology, Haifa 32000, Israel}, \and Paul Anthony Smith\footnote{UCLA Mathematics Department, Los Angeles, CA  90095-1555, USA}}

\maketitle

\pagebreak

\begin{abstract}
Characterized are all simple undirected graphs $G$ such that any
real symmetric matrix that has graph $G$ has no eigenvalues of
multiplicity more than 2.  All such graphs are partial 2-trees
(and this follows from a result for rather general fields), but only certain
partial 2-trees guarantee maximum multiplicity 2. Among partial
linear 2-trees, they are only those whose vertices can be covered
by two "parallel" induced paths. The remaining graphs that
guarantee maximum multiplicity 2 are comprised by certain
identified families of "exceptional" partial 2-trees that are not
linear.
\end{abstract}


\nline


\pagebreak

\section{Introduction}
\setcounter{equation}{0}

Throughout, $G$ denotes a simple, undirected graph on $n$ vertices
without loops.  Associate with $G$ the set $S(G)$ of all $n$-by-$n$,
real symmetric matrices $A$ whose graph is $G$.  No restriction
(other than reality) is placed upon the diagonal entries of $A$ by
$G$.  For each $A \in S(G)$, let $M(A)$ be the largest multiplicity
for an eigenvalue of $A$ and let rank($A$) denote the rank of $A$.
Then, over $A \in S(G)$,
$$M(G) = \text{max }M(A)$$
and
$${\rm msr}(G) = \text{min rank }A,$$
the maximum multiplicity for
$G$ and the minimum symmetric rank in $S(G)$, respectively.
Because all eigenvalues of matrices in $S(G)$ are real and translation
by a real multiple of the identity does not change membership in
$S(G)$, of course
$$M(G) + {\rm msr}(G) = n,$$
and the two may be viewed interchangeably. This allows us to
implicitly assume, that, when working with $M(A)$ and when
convenient, 0 is an eigenvalue that attains $M(A)$. The same holds
for $M(G)$. When $\mathbb{R}$ is replaced by a field F, then
$M(G)$ is defined as the maximum corank of all symmetric matrices
with entries in F, and whose graph is G.

In \cite{F}, it was observed that the only graph $G$ for which
$M(G) = 1$ is the path on $n$ vertices. In \cite{JLD}, the maximum
multiplicity $M(G)$ has been determined whenever $G$ is a tree.
Our purpose here is to describe all graphs $G$ for which $M(G) =
2$, a much larger (than $M(G) = 1$) and more subtle to describe
class.

\section{Partial 2-trees and preliminaries}
\setcounter{equation}{0}

Recall that a {\em k-tree} is a graph sequentially constructed
from $k+1$-cliques ($K_{k+1}$) via articulation along $k$-cliques.
Thus, a traditional tree is a 1-tree.  We are particularly
interested here in 2-trees, in which the building blocks are
triangles ($K_3$'s) and the articulation is along edges.  A {\em
partial $k$-tree} is a $k$-tree from which some edges (without
incident vertices) have been deleted.  We call a 2-tree linear if
it has precisely two vertices of degree two; we also consider
$K_3$ to be a (degenerate) linear 2-tree. In this event, there is
a natural order to the triangles and a linear 2-tree is somewhat
analogous to a path (though it should be noted that a linear
2-tree may have vertices of arbitrarily high degree).

A graph $H$ is a {\em homeomorph} of a graph $G$ if $H$ may be
obtained from $G$ by a sequence of edge subdivisions.  We use
$hK_4$ and $hK_{2,3}$ to denote graphs that are homeomorphs of
$K_4$ and $K_{2,3}$ (the complete bipartite graph on two and three
vertices) respectively.  An $hK_{2,3}$ is just the result of
articulation of two cycles along a common induced path of at least
two edges.

\noindent{\bf Examples.} Let
$$
A_1 =
\begin{bmatrix}
1 & 1 & 1 & 1\\
1 & 1 & 1 & 1\\
1 & 1 & 1 & 1\\
1 & 1 & 1 & 1
\end{bmatrix} \in S(K_4)
$$
and
$$
A_2 =
\begin{bmatrix}
-1 & 0 & 1 & 1 & 1\\
0 & 1 & 1 & 1 & 1\\
1 & 1 & 0 & 0 & 0\\
1 & 1 & 0 & 0 & 0\\
1 & 1 & 0 & 0 & 0
\end{bmatrix} \in S(K_{2,3}).
$$
Then ${\rm rank}(A_1) = 1$ and ${\rm rank}(A_2) = 2$, and it is easy
to see that $M(K_4) = M(K_{2,3}) = 3$.  Note that the Schur
complement (see \cite{HJ1}, Ch.~0) in $A_2$ of the $(1,1)$ entry gives $A_1$ and, thus,
something whose graph is $K_4$.

\begin{La}\label{L2-1}
Let $G^\prime$ be the graph resulting from an edge subdivision in
the graph $G$.  Then $M(G^\prime) = M(G)$ or $M(G^\prime) = M(G) +
1$ (i.e., $M(G^\prime) \geq M(G)$).

\end{La}
\begin{proof} Denote by $e = (v_1,v_2)$ the edge in $G$ that is
subdivided to obtain
$G^\prime$.  After subdividing $e$, we get a new vertex $v$
whose only neighbors are $v_1$ and $v_2$.
Let us number vertices $v_1$, $v_2$, and $v$ by the numbers
$n-1$, $n$, and $n+1$, respectively.
Here and in the sequel we shall assume that if some vertices
of a graph $G$ have been numbered, then any matrix in
$S(G)$ that we consider is consistent with the numbering
(we shall only use integers in the set $\{1, \ldots, n\}$).
Note that by permutation similarity we may always transform
an arbitrary matrix $B \in S(G)$ to one consistent with a
numbering.
Let $A \in S(G^\prime)$
satisfy $M(A) = M(G^\prime)$, i.e.,
${\rm rank}\text{ } A = (n+1) - M(G^\prime)$.
We split the proof into two (mutually exclusive) cases:
\begin{eqnarray*}
&(a)& \text{ the $(n+1)^{st}$ diagonal entry of $A$ is nonzero.} \\
&(b)& \text{ the $(n+1)^{st}$ diagonal entry of $A$ is zero.}
\end{eqnarray*}
Let us first suppose that our $A$ as defined above satisfies
condition $(a)$. Only the last two off-diagonal entries of the
$(n+1)^{st}$ row and of the $(n+1)^{st}$ column are nonzero. We
may therefore add multiples of the $(n+1)^{st}$ column of $A$ to
columns $n-1$ and $n$ so that the entry in the last row of each
column is zero. By symmetry we may simultaneously perform the same
operation with the roles of rows and columns interchanged.  Call
the matrix we so obtain $\tilde{A}$. As a result of our
operations, $\tilde{A}$ is a direct sum of a (real symmetric)
matrix $B$ with graph $G$ and a single nonzero number $x$, i.e.,
\begin{equation*}
\tilde{A} = \begin{bmatrix}B & 0 \\ 0 & x\end{bmatrix}.
\end{equation*}
Since $A$ was chosen to be of minimum
possible rank, it follows that $B$ has minimum possible rank also, and
so ${\rm rank}\text{ } B = n - M(G)$.  Therefore
the chain of equalities
\begin{equation}
n - M(G) + 1 = {\rm rank }\text{ }B + 1 = {\rm rank}\text{ }A = (n+1) - M(G^\prime)
\label{2.1}
\end{equation}
holds, whence $M(G^\prime) = M(G)$.

Now let us suppose that $A$ satisfies case $(b)$
and moreover that
there is no matrix $A^\prime \in S(G)$ with $M(A^\prime) = M(G)$ that
satisfies condition $(a)$.
Add $1$ to the $(n+1)^{st}$ diagonal entry of $A$
and call the new matrix $\tilde{A}$.
Due to our assumption, ${\rm rank}\text{ } \tilde{A} =
1 + {\rm rank}\text{ } A$, and so
\begin{equation}
{\rm rank}\text{ } \tilde{A} = (n+1) - M(G^\prime) + 1.
\label{2.2}
\end{equation}
We now apply to the matrix $\tilde{A}$ the procedure
used to prove part $(a)$.
By (\ref{2.1}) and (\ref{2.2}), we obtain
$M(G^\prime) = M(G) + 1$.
\end{proof}

\noindent{\bf Remark.} We note that both eventualities may occur.
If $G$ is a cycle, then $M(G^\prime) = M(G)$, and if
$G$ consists of two cycles that overlap in one (and only one) edge,

\includegraphics[width=70mm]{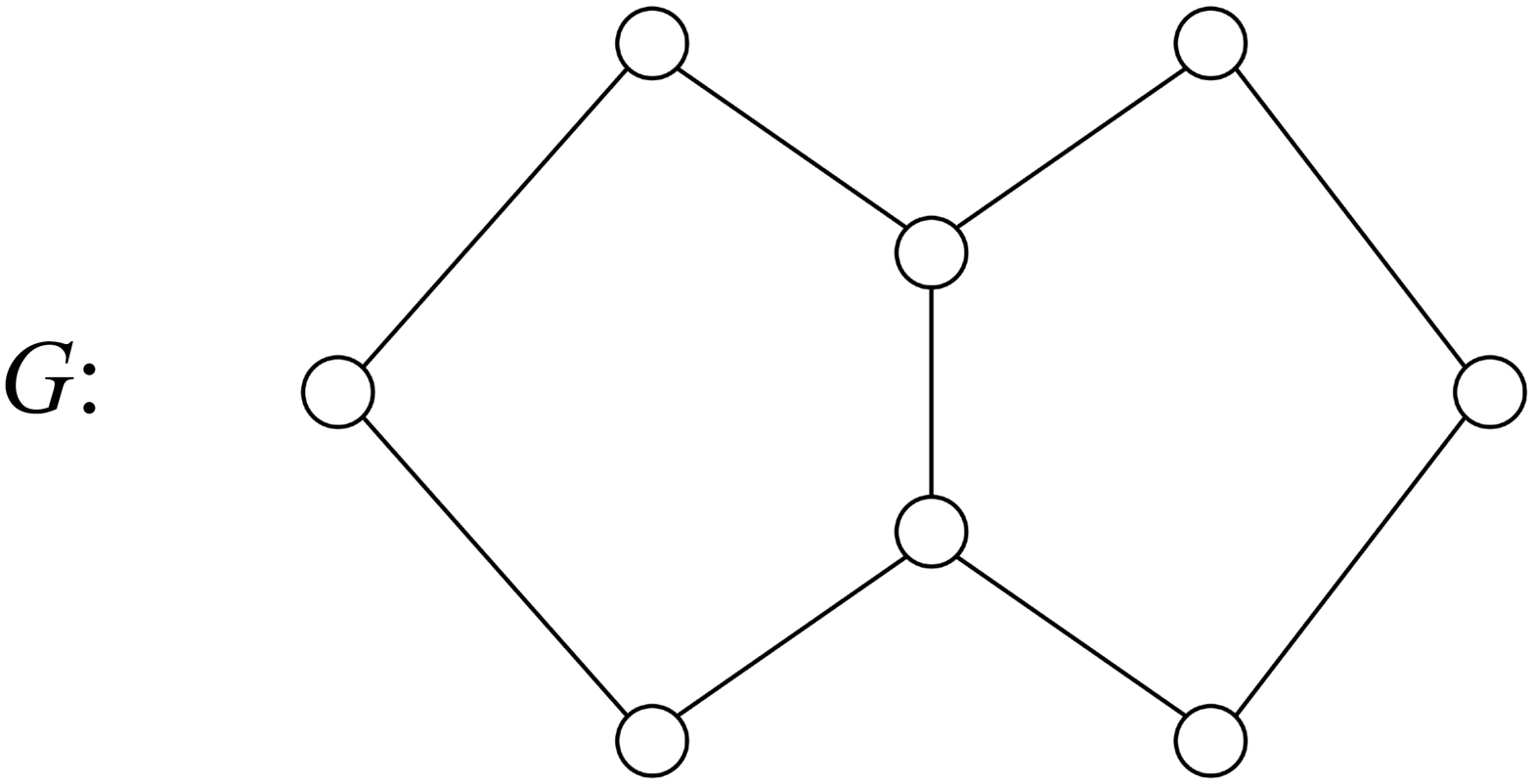}

\noindent then
$M(G^\prime) = M(G) + 1$ if the overlapping edge is subdivided.

\vspace{\baselineskip}

Because of Lemma \ref{L2-1}, we see that any graph $G$ that is either an
$hK_4$ or an $hK_{2,3}$ satisfies $M(G) \geq 3$.

The following combinatorial characterization of partial 2-trees
is known (see \cite{WC} or \cite{BLS})
 and will be useful to us.

\begin{La}
The graph $G$ is a partial 2-tree if and only if $G$ does not
contain an induced subgraph that is a supergraph of an $hK_4$.
\label{L2-2}
\end{La}

We may now establish a key step in our characterization of graphs
for which $M(G) = 2$. Our proof of a statement for more general
fields is given in Appendix A.

\begin{La}
If $G$ is a graph such that $M(G) = 2$, then $G$ is a partial
2-tree.
\label{L2-3}
\end{La}

In fact, a stronger result than Lemma \ref{L2-3} holds (cf.
\cite{vdH} ), namely: if G is not a partial 2-tree, then there
exists a positive semi-definite matrix in $S(G)$ with nullity
$\geq 3$. But this result has no natural analogue over general
fields.

Of course, not all partial 2-trees have maximum multiplicity two.
For example, $K_{2,3}$ is a partial 2-tree (simply add an edge
between the two vertices in the first part to produce a ``book''
of triangles articulated at a single edge, a graph for which the
maximum multiplicity is also greater than two). The rest of our
work is to sort out which partial 2-trees do have maximum
multiplicity two. In the next section, we identify the major
portion of them, but certain ``exceptions'' will be identified
later.

\section{Graphs of two parallel paths}
\setcounter{equation}{0}

\begin{Dn}
A graph $G$ is a graph of two parallel paths if there exist two
independent induced paths of $G$ that cover all the vertices of
$G$ and such that any edges between the two paths can be drawn
so as to not cross.  A simple path is not considered to be
such a graph (and two paths not connected is considered to be
such a graph).
\end{Dn}
We shall call two independent induced paths satisfying the conditions
in the above definition a {\em pair of parallel paths}.
We note that $K_3$ is a graph of two parallel paths, and,
in any given pair of two parallel paths of $K_3$, one of these paths
is degenerate (a vertex).   We note also that each graph of two
parallel paths is a partial linear 2-tree.  This we shall prove
after we elaborate on the meaning of the
the requirements in the definition.

\begin{Rk}
The matrix structure of graphs of two parallel paths: \label{R3-2}
\end{Rk}

Here we express more precisely the requirement in the above
definition that a pair of parallel path may be drawn so that edges
between the path do not cross. Suppose we have a graph $G$ on $n$
vertices such that there exists a pair of independent paths $P_1$
and $P_2$. Let $k_i$ denote the number of vertices of $P_i$.
Number the vertices of $P_1$ consecutively from $1$ to $k_1$,
starting from a pendant vertex of $P_1$.  We shall number the
vertices of $P_2$ similarly, but in this case we shall require
that we start numbering from a pendant vertex of $P_2$ such that
we minimize the number of times the following situation occurs: a
vertex $j$ of $P_2$ is adjacent to a neighbor $s$ of $P_1$ and
some vertex $k
> j$ of $P_2$ is adjacent to some $t < s$ of $P_1$.  The vertices
of $P_2$ may be numbered so that this situation never occurs if
and only if $G$ is a graph of two parallel paths.  This may be
rephrased as following. $G$ is a graph of two parallel paths if
and only if there exists $A \in S(G)$ of the following form:
\begin{equation*}
A = \begin{bmatrix}T_1 & B \\ B^t & T_2\end{bmatrix},
\end{equation*}
where $T_1$ and $T_2$ are irreducible and tridiagonal and $B$
satisfies the following:
\begin{equation*}
\begin{array}{l}
\text{Whenever $b_{ij} \neq 0$ for some entry, then $b_{kl} = 0$}\\
\text{for $k > i$ and $l < j$, and for $k < i$ and $l > j$}.
\end{array}
\end{equation*}
and in addition $B$ is such that whenever $B \neq 0$ and
$b_{k_1,k_1+1} \neq 0$, this entry is not
the only nonzero entry of $B$ (this excludes
paths).

\begin{Rk}
If $G$ is a graph of two parallel paths, and $P_1$ and $P_2$
constitute a pair of two parallel paths, then we shall assume that
$P_1$ and $P_2$ are numbered so that vertices of $P_1$ have lower
numbers than those of $P_2$, and so that $A \in S(G)$ is of the
form given in Remark \ref{R3-2}.
\end{Rk}

\begin{La}
If a graph $G$ is a graph of two parallel paths, then $G$ is a partial
linear 2-tree.
\end{La}
\begin{proof}
Let $G$ be a graph of two parallel paths on $n$ vertices, and let
$p_1$ and $p_2$ be a pair of parallel paths, with $n_1$ and $n_2$
vertices, respectively.  We shall construct a finite sequence of
graphs, the last of which we shall show must be a linear 2-tree.
Suppose $H$ is a graph of two parallel paths with a pair of paths
$q_1$ and $q_2$.  Let $v$ be a vertex in $q_1$ numbered with $j$.
We shall define a function ${\rm upper}_H(j)$ on the vertices of
$q_1$ as follows.  Let $S$ be the set of vertices of $q_1$ with
label $j+1$ or greater.  Let $S^\prime$ be the subset of all
vertices in $S$ that are adjacent to some vertex in $q_2$. If
$S^\prime$ is empty, define ${\rm upper}_H(j)$ to be the last
vertex in $q_2$.  If $S^\prime$ is not empty, take the vertex in
this set with the lowest number, call it $w$, and define ${\rm
upper}_H(j)$ to be the lowest numbered neighbor of $w$.  Define
${\rm lower}_H(j)$ similarly, but with the roles of ``greatest''
and ``lowest'' interchanged; S is now defined as the set of
verices of $q_1$ with label $j-1$ or lower, and $S^\prime$ is
defined in the obvious way. Let us now construct the following
sequence of graphs. Let $G_0 = G$ and, for $i = 1, \ldots, n_1$,
let $G_i$ be the supergraph of $G_{i-1}$ obtained by articulating
edges $(i,j)$, where $j$ runs from ${\rm lower}_{G_{i-1}}(i)$ to
${\rm upper}_{G_{i-1}}(i)$. It may happen that some of these edges
already exist.  Note that the graph $G_{n_1}$, by construction,
consists of only triangles, and so it a 2-tree.  Moreover, either
$1$ or $n_1 + 1$ is of degree two, and either $n_1$ or $n$ is of
degree two, and all other vertices are of degree at least three.
To see this, note that not both $1$ and $n_1 + 1$ can be of degree
greater than two, for in that case $p_1$ and $p_2$ could not be
drawn so that edges do not overlap. If, in $G$ we do not already
have that each of these vertices is of degree at least two, then,
by considering the remaining cases, it is easy to see that our
construction gives the result claimed. We may argue similarly for
vertices $n_1$ and $n$. On the other hand, any other vertex in
$p_1$ must be of degree at least three, since ${\rm
lower}_{G_{i-1}}(j) \leq {\rm upper}_{G_{i-1}}(j)$ for each $i$
and $j$.  It is also easy to see that each vertex of $p_2$ of
degree two relative to $p_2$ is numbered such that, for some $i$,
$j$, ${\rm lower}_{G_{i-1}}(j) \leq p_2 \leq {\rm
upper}_{G_{i-1}}(j)$. Therefore $G_{n_1}$ is a 2-tree with
precisely two vertices of degree two, whence a linear a 2-tree.
\end{proof}

\begin{Dn}
We say that a graph is $C_2$ if it is connected and no vertex has
degree less than two (no pendant vertices).
\end{Dn}

\begin{Dn}
The core of a connected graph $G$ is the maximal induced subgraph that is
$C_2$.
\end{Dn}

\begin{La}
Suppose that $G$ is a graph of two parallel paths on $n$ vertices.
Then ${\rm msr}(G) = n - 2$, or, equivalently, $M(G) = 2$.
\label{L3-7}
\end{La}
\begin{proof}
Suppose that the two parallel paths are $P_1$ and $P_2$ with $k_1$
and $k_2$ vertices respectively, $k_1 + k_2 = n$.  By definition,
$k_1,k_2 \geq 1$.  We number the vertices of $P_1:1,2,\ldots,k_1$
consecutively along the path and the vertices of $P_2:k_1 + 1,
\ldots,n$ also consecutively along the path, but beginning with
$k_1 + 1$ at the same end of $P_2$ as $1$ is of $P_1$ (if
unambiguous, and at either end otherwise).  Let $A \in S(G)$. We
show by induction on $n$ that $B=A(\{k_1,n\},\{1,k_1+1\})$, the
submatrix of $A$ obtained by deleting rows $k_1$ and n and columns
1 and $k_1+1$, is permutation equivalent to a triangular matrix
with nonzero diagonal.  In the event that $k_1 = 1$ or $k_2 = 1$,
this is immediate, as then the other is $n - 1$ and the indicated
submatrix is a necessarily triangular submatrix whose diagonal is
nonzero because of an irreducible tridiagonal principal submatrix.
Thus, we may assume that $k_1,k_2 \geq 2$ and that the cases $n =
2,3$ have been verified.  To start the induction, the case $n = 4,
k_1 = 2 = k_2$, is also easily checked.  Because of the ``no
crossing'' requirement upon edges between $P_1$ and $P_2$, not
both vertices $1$ and $k_1 + 1$ can have degree more than $2$ and
not both vertices $k_1$ and $n$ can have degree more than $2$.
Thus, in the matrix $B$, either column $k_1$ or $n$ (original
numbering) has exactly one nonzero entry. We may assume without
loss of generality that it is column $n$ and that it appears in
the last (current numbering) position of that column of $B$.
Deletion of this row from $B$ leaves, either a path, or, by the
induction hypothesis, a matrix containing an $(n-3)-$triangle,
which is extended to an $(n-2)-$triangle by the nonzero in the
last column of $B$, to complete the induction.

We have shown that msr$(G)\geq n-2$. In fact, the method of proof,
which is purely combinatorial, shows that the same conclusion
holds for any field F.

To complete the proof, we note that G is not a path and apply
\cite{F} (or, for any infinite field, apply \cite{BD}) to conclude
that msr$(G)\leq n-2$.
\end{proof}

 \noindent{\bf Example.}
 Let $G$ be the unique 5-vertex linear 2-tree.  Let $P_1$
and $P_2$ be two parallel paths and suppose that $P_1$
has two vertices.  Furthermore, suppose that the first
and last vertices of $P_2$ are each of degree two.
In this example we shall use
$P(A^\prime)$ to denote the pattern of $A^\prime$.
Numbering according to our prescription
implies that matrices in $S(G)$ consistent with this numbering
have the pattern
\begin{equation*}
P(S(G)) :=
\begin{bmatrix}
\cdot & \times & \times & \times & 0 \\
\times & \cdot & 0 & \times & \times\\
\times & 0 & \cdot & \times & 0\\
\times & \times & \times & \cdot & \times\\
0 & \times & 0 & \times & \cdot
\end{bmatrix}
\end{equation*}
where $\times$ denotes a nonzero entry and $\cdot$ a completely
free entry.  Let $A \in S(G)$ satisfy $M(A) = M(G)$.  In
particular, $A$ has the sign pattern given above, i.e. P(A) =
P(S(G)). The pattern of  $P_1$ is given by the upper diagonal $2
\times 2$ block in $P(A)$, and $P_2$ by the lower diagonal $3
\times 3$ block. The pattern of $B$ is given by the uppermost
rightmost $2 \times 3$ block in $P(G)$. Striking out rows 2 and 5
and columns 1 and 3, we obtain a submatrix $\tilde{A}$ with
pattern
\begin{equation*}
P(\tilde{A}) =
\begin{bmatrix}
\times & \times & 0\\
0 & \times & 0\\
\times & \cdot & \times
\end{bmatrix}
\end{equation*}
By taking the Schur complement with respect to the (3,3) nonzero
entry, we see that $\tilde{A}$ has full rank, and therefore
${\rm msr}(A) \geq 3$, or equivalently, $M(A) < 3$.  Since $A$ is not
a path, $M(A) > 1$, and therefore we conclude that $M(A) = 2$.

\section{Graphs of minimum degree two and M = 2}
\setcounter{equation}{0}

\begin{La}
If a $C_2$ graph $G$ contains a cut-vertex, then $M(G) > 2$.
\label{L4-1}
\end{La}
\begin{proof} Let us denote the cut-vertex by $v$ and consider
the induced subgraph $G- v$.  We are left with connected components
$K_1, \ldots, K_n$, where $n$ is at least two.
Let us introduce the induced subgraphs
\begin{equation*}
\tilde{K}_1 = G-\{K_2, \ldots, K_n\}
\end{equation*}
and
\begin{equation*}
\tilde{K}_2 = G-K_1.
\end{equation*}
Neither $\tilde{K}_1$ nor $\tilde{K}_2$ is a path.  By \cite{F} ,
\cite{BD}, the lemma holds for every infinite field $F$.
\end{proof}

The following lemma is a special case of Theorem 2.3 in \cite
{BFH}, and is brought here for the sake of completeness. It holds
for any field F.
\begin{La}
Let $G$ be a graph containing a pendant vertex $v$ with unique
neighbor $u$.  Then
$M(G) = {\rm max}\left\{M\left(G-v\right), M\left(G -
    \{u,v\}\right)\right\}$.
\label{L4-2}
\end{La}
\begin{proof} Number the vertex $v$ by $n$ and $u$ by $n-1$.
Let $A \in S(G)$ satisfy $M(A) = M(G)$.  We consider
separately two
(mutually exclusive) cases:
\begin{eqnarray*}
&(a)& \text{ the $n^{th}$ diagonal entry of $A$ is nonzero.} \\
&(b)& \text{ the $n^{th}$ diagonal entry of $A$ is zero.}
\end{eqnarray*}
Let us suppose that we may find a matrix $A$ as defined above
such that condition $(a)$ holds.  Taking into account that
the only nonzero entries in the last row are entries $n-1$ and $n$,
we add multiples of the last row of $A$ to the row $n-1$
so that its $n^{th}$ entry becomes zero.
We simultaneously perform the same procedure with the roles
of rows and columns reversed.  Let us call the resulting matrix
$\tilde{A}$.  As a result
of our operations, $\tilde{A}$ is a direct sum of a
(real symmetric) matrix $B$ with graph $G-v$ and a single
nonzero number $x$:
\begin{equation}
\tilde{A} = \begin{bmatrix}B & 0\\0 & x\end{bmatrix}.
\label{4.1}
\end{equation}
Note that $M(A) = M(G)$ forces $M(B) = M(G-v)$: clearly $M(B) \leq
M(G - v)$; if $M(B) < M(G - v)$, then ${\rm msr}(G) > {\rm msr}(G
- v) + 1$, a contradiction. Because $x \neq 0$ in our direct sum
decomposition, it holds that $M(G) = M(A) = M(B) = M(G-v)$.

Now suppose that $A$ satisfies condition $(b)$ and that there
is no matrix $A^\prime \in S(G)$ with $M(A^\prime) = M(G)$
that satisfies condition $(a)$.
Since $(b)$ is satisfied, in the last row of $A$ only the $(n-1)^{st}$
entry is nonzero.  We therefore  may add multiples of the last
row to each row of $A$, canceling all nonzero entries
of the $(n-1)^{st}$ column of $A$ without affecting any other
entries.  By symmetry, we may perform the same operation
with the roles of rows and columns reversed.  Let us
call the resulting matrix $\hat{A}$.  Note that $\hat{A}$
is a direct sum of a (real symmetric) matrix $B$ with graph
$G-\{u,v\}$ and a $2\times 2$ matrix $X$ given by
$X = \begin{bmatrix}0 & x \\ x & 0\end{bmatrix}$, where $x$
is nonzero:
\begin{equation}
\hat{A} = \begin{bmatrix}B & 0 \\ 0 & X\end{bmatrix}.
\label{4.2}
\end{equation}
Since $M(A) = M(G)$, it follows that $M(B) = M(G-\{u,v\})$.
By our direct sum decomposition and the fact that $X$ has
full rank, we get that $M(G) = M(A) = M(B) = M(G-\{u,v\})$.
Thus we have established that
$M(G) \in \left\{M(G-v), M(G-\{u,v\}\right\}.$  Now, suppose
that we start with a matrix $\tilde{A}$ as given in
(\ref{4.1}) where $B \in S(G-v)$, $M(B) = M(G-v)$,
and $x$ is nonzero.  By reversing all of our row and column
operations, we may obtain a matrix $A \in S(G)$ with $M(A) = M(G-v)$.
Similarly, if we start with a matrix $\hat{A}$ as given
in (\ref{4.2}) where $B \in S(G-\{u,v\})$,
$M(B) = M(G-\{u,v\})$, and $X = \begin{bmatrix}0 & x\\x & 0
\end{bmatrix}$, where $x \neq 0$, then by reversing our row
and column operations there performed, we may obtain a matrix
$A \in S(G)$ with $M(A) = M(G-\{u,v\})$.
\end{proof}

\begin{La}
Let $G$ be a graph containing an induced subgraph that is a
supergraph of an $hK_{2,3}$.  Then $M(G) \geq 3$.
\label{L4-3}
\end{La}
\begin{proof}
Let the $hK_{2,3}$ consist of three internally independent paths,
of at least two edges each between vertices $u$ and $v$. Call them
$P_1$, $P_2$, and $P_3$.  If $G$ contains an $hK_4$, then the
conclusion follows from Lemma \ref{L2-3}.  Thus, we may assume
that there is no path in $G-\{u,v\}$ from an interior vertex of
$P_i$ to an interior vertex of $P_j$, $j \neq i$, $1 \leq i,j \leq
3$. It follows that if $u$ and $v$ are deleted from $G$, then at
least three components result, including a component
``corresponding'' to each $P_i$, $i=1,2,3$.  We conclude that $A
\in S(G)$ appears, with proper numbering of vertices, as
\begin{equation*}
A =
\begin{bmatrix}
\cdot & ? & d_1 & d_2 & d_3 & ? \\
? & \cdot & f_1 & f_2 & f_3 & ? \\
d_1^t & f_1^t & A_1 & 0 & 0 & 0\\
d_2^t & f_2^t & 0 & A_2 & 0 & 0\\
d_3^t & f_3^t & 0 & 0 & A_3 & 0\\
? & ? & 0 & 0 & 0 & B
\end{bmatrix}
\end{equation*}
Here $A_i$ corresponds to component $G_i$ and includes the interior
vertices of $P_i$, $i=1,2,3$, $B$ may be empty, and the first two
rows and columns correspond to $u$ and $v$.  Each $d_i$ has its
first entry nonzero, and each $f_i$ has its last
entry nonzero.  The symbol ``?'' is either free nonzero or zero.
Now, identify the three
neighbors of $u$ along the path $P_1$, $P_2$, and $P_3$ together with
$u$ and $v$ to give the set $S$ contained among the vertices of $G$.
Let $A \in S(G)$ be chosen as follows: $A[S^c]$ is an $M$-matrix
(see \cite{HJ2}, Ch.~2),
the three edges from $u$ to its neighbors in $S$ are free, the
diagonal entries corresponding to $S$ are free and all other edges
in $G$ are chosen to be positive.  The Schur complement
$A / A[S^c]$ then appears as
\begin{equation*}
\left[
\begin{array}{c|c}
\begin{matrix}\tilde{a}_{11} & ? \\ ? & \tilde{a}_{22}\end{matrix} &
\begin{matrix}\tilde{a}_{13} & \tilde{a}_{14} & \tilde{a}_{15} \\ a & b & c\end{matrix}\\
\hline
\begin{matrix}\tilde{a}_{13} & a \\ \tilde{a}_{14} & b \\ \tilde{a}_{15} & c\end{matrix} &
\begin{matrix}\tilde{a}_{33} & 0 & 0 \\ 0 & \tilde{a}_{44} & 0 \\
0 & 0 & \tilde{a}_{55} \end{matrix}
\end{array}
\right]
\end{equation*}
in which ? is either free nonzero, fixed nonzero, or zero, $\tilde{a}_{13}$, $\tilde{a}_{14}$,
$\tilde{a}_{15}$ are free nonzero, $a$, $b$, $c$ are nonzero and
$\tilde{a}_{11}, \ldots, \tilde{a}_{55}$ are free.

To see that $a$, $b$, and $c$ may be chosen to be nonzero, consider
the following two situations.  Denote by $w$ the neighbor of $u$
along $P_1$.  If $w$ is not a neighbor of $v$, then the sign of
the second summand in the Schur complement $A / A[S^c]$ shows that
$a < 0$.  If $w$ is a neighbor of $v$, then $a$ is given
as the difference of two positive numbers.  In this case we may
treat the entry corresponding to the edge $(w,v)$ as a free variable,
and therefore guarantee that $a \neq 0$.  The same argument holds for
$b$ and $c$.

If we choose
$\tilde{a}_{33} = \tilde{a}_{44} = \tilde{a}_{55} = 0$, $\begin{bmatrix}\tilde{a}_{13}
& \tilde{a}_{14} & \tilde{a}_{15}\end{bmatrix}$
proportional to $\begin{bmatrix}a & b & c\end{bmatrix}$ and either
$\begin{bmatrix}\tilde{a}_{11} & ? \\ ? & \tilde{a}_{22}\end{bmatrix}$ to be
0 if $? = 0$ or to be so that
$\begin{bmatrix}\tilde{a}_{11} & ? & \tilde{a}_{13} \\
? & \tilde{a}_{22} & a\end{bmatrix}$ is rank one, otherwise, the proof
is complete.
Note that, because of the Schur complement step, $\tilde{a}_{13}$,
$\tilde{a}_{14}$, $\tilde{a}_{15}$ and $?$ could each have a single forbidden
value other than zero.  This does not debilitate the argument
because of the flexibility in proportionality and in the
upper left 2-by-2 block.
\end{proof}
\begin{Rk}
Lemma \ref{L4-3} holds for any infinite field. Indeed, a
straightforward computation shows that the entries in the
$3,4;3,5$ and $4,5$ positions of $A/A[S^c]$, as well as the
entries in the transpose positions are 0. The results about the
other entries of $A/A[S^c]$ follow the discussion in Appendix A,
and in particular Lemma \ref{LA-2}
\end{Rk}
\begin{Dn}
A graph $G$ is a $SEAC$ (singly edge-articulated cycle) graph if
it is sequentially built from cycles via articulation of each
successive cycle along an edge of the previous graph.  No such edge
may be used more than once for articulation.
\end{Dn}

\begin{La}
If $G$ is a $C_2$ graph such that $M(G) = 2$, then $G$ is a $SEAC$
graph.
\label{L4-5}
\end{La}
\begin{proof}
By Lemmas \ref{L2-3} and \ref{L4-1} we can assume that $G$ is a
partial 2-tree which does not contain any cut vertex. The lemma is
proved by induction on $n=|G|$. It holds for $n=3$, so we consider
the general induction step.\\

G is a partial 2-tree, so it has a supergraph $\tilde{G}$ with
$|\tilde{G}|=|G|=n$, and such that $G$ is obtained from
$\tilde{G}$ by removing some edges. $\tilde{G}$ consists of $n-2$
triangles. Denote by $T_{n-2}$ the last one that was articulated
in the construction of $\tilde{G}$. Denote its vertices by
$u,v,w$, where we may assume that the degree of w in $\tilde{G}$
is 2.\\ \\
\centerline{\includegraphics[width=70mm]{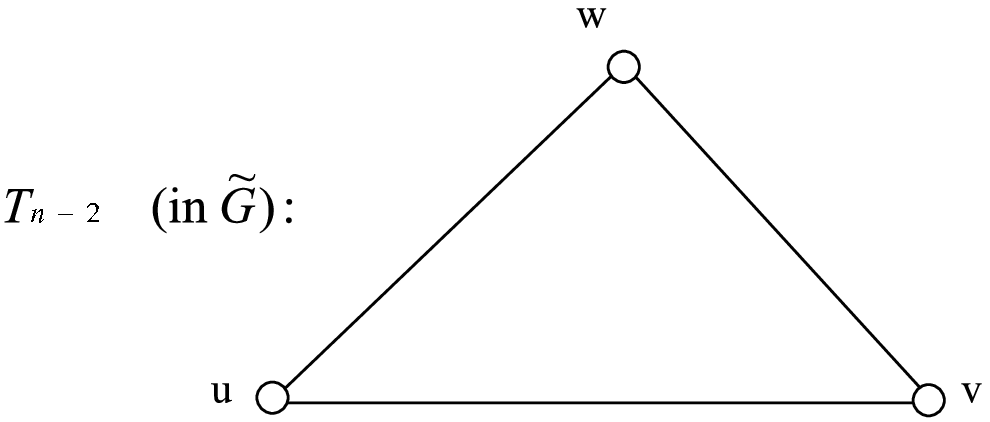}}\\

Since $G$ contains no cut vertices, $uw$ and $vw$ are edges of
$G$.%
\Cs Suppose that $uv$ is an edge of $G$. Then, since $G$ contains
no cut vertices, the degrees of $u$ and $v$ in $G$ are at least 3.

Let $G'=G-\{w\} .$ Clearly, $G'$ is a $C_2$ graph with $|G'|=n-1$.
By \cite{F} we cannot have $\text{msr}(G')=(n-1)-1 .$ We cannot
have $\text{msr}(G')\leq (n-1)-3$, because this implies
$\text{msr}(G)\leq n-3$. We use here $\text{msr}(K_3)=1$. Hence
$\text{msr}(G')=n-3=(n-1)-2$, so we can apply the induction
hypothesis and conclude that $G'$ is a SEAC. We are done if $uv$
is incident to exactly one cycle in $G'$, and this is indeed the
case by Lemma \ref{L4-3}.%
\Cs Suppose that $uv$ is not an edge of $G$. Let $G'$ be obtained
from $G$ by compressing $w$, so $|G'|=n-1$, and $uv$ is an edge of
$G'$. $G'$ is $C_2$, and as in Case 1, msr$(G')=(n-1)-2$, so the
induction hypothesis implies that $G'$ is a SEAC. If the edge $uv$
(in $G'$) is incident to only one cycle we are done. Otherwise,
going back to $G$ we get a contradiction by Lemma \ref{L4-3}.
\end{proof}

\begin{Dn}
In a $SEAC$ graph, the cycles used to build it are well defined, as
are the edges of articulation.  Each edge of articulation uniquely
defines two of the cycles.  We say that two of the cycles are
neighbors if they share an edge of articulation.  A $SEAC$ graph is
called linear $(LSEAC)$ if each of its constituent cycles has at
most two neighbors.  An $LSEAC$ graph that consists of more than a
cycle has just two cycles with only one neighbor each (the two ends
of the linear path of cycles).
\end{Dn}

\begin{La}
If $G$ is a $C_2$ graph and $M(G) = 2$, then $G$ is an $LSEAC$ graph.
\label{L4-7}
\end{La}
\begin{proof} Let $G$ be a $C_2$ graph with $M(G) = 2$.
By Lemma \ref{L4-5}, $G$ is a $SEAC$ graph.
Thus it remains to be shown that $G$ is in
fact an $LSEAC$ graph.  If $G$ is not an $LSEAC$ graph, then
there exists a cycle $Z$ with at least $k \geq 3$ neighbors, say
$Z_1, \ldots, Z_k$ ($k \geq 3$).
These neighbors uniquely determine
the connected components of the graph $G - Z$.  Let us superimpose
the connected component corresponding to $Z_i$
with the intersection of $Z_i$ and $Z$, and denote it by
$\tilde{Z}_i$.  It is clear $\tilde{Z}_i$ is a $SEAC$ graph.
We shall assume that each is in fact an $LSEAC$ graph.  Once
the lemma is proved for this case, the general case follows
by induction.  Let $s$ denote the number of vertices of $Z$
and $s_i$ the number of vertices of $\tilde{Z}_i$.  The cycle
$Z$ contains $s-k$ edges not shared with some $\tilde{Z}_i$.
It therefore follows that
\begin{equation}
{\rm msr}(G) \leq {\rm msr}(\tilde{Z}_1) + \ldots + {\rm msr}(\tilde{Z}_k) + (s-k).
\label{4.3}
\end{equation}
Since each $\tilde{Z}_i$ is an $LSEAC$ graph,
\begin{eqnarray}
{\rm msr}(\tilde{Z}_1) + \ldots + {\rm msr}(\tilde{Z}_k) + (s-k) &=&
 (s_1 - 2) + \ldots + (s_k - 2) + (s-k) \nonumber \\
&=& (n-s) + (s-k) = n - k,
\label{4.4}
\end{eqnarray}
By (\ref{4.3}) and (\ref{4.4}), ${\rm msr}(G) \leq n - k$, i.e.,
$M(G) \geq k$.  In particular, $M(G) > 2$.
The general case follows by induction.
\end{proof}

\begin{Tm}
If $G$ is a $C_2$ graph, then the following three statements are
equivalent:
\begin{enumerate}
\item M(G) = 2
\item G \text{ is a graph of two parallel paths, and }
\item G \text{ is an LSEAC graph.}
\end{enumerate}
\label{T4-8}
\end{Tm}
\begin{proof} Note that $2. \implies 1.$ follows from Lemma
\ref{L3-7} and that $1. \implies 3.$ follows from Lemma
\ref{L4-7}.  It remains to be shown that $3. \implies 2.$.
This is trivial if $G$ is a single cycle, and so let us
suppose that $G$ contains at least two cycles.  Let
$Z_s$ and $Z_t$ denote the two cycles with only one neighbor
each.  Let $s_1, s_2$ be two adjacent vertices belonging
to $Z_s$ such that the edge that they determine is not an
edge of articulation.  Define $t_1, t_2$ belonging to $Z_t$
similarly.
Now remove from $G$ all edges of articulation.  The resulting
graph is a cycle that clearly defines the two
parallel paths in $G$.
\end{proof}

\section{The exceptional graphs and the general result}
\setcounter{equation}{0}

In this section we shall prove the following general result:

\begin{Tm}
The graph $G$ satisfies $M(G) = 2$ if and only if $G$ is a graph
of two parallel paths or $G$
is one of the types listed in Table B1.
\label{T5-1}
\end{Tm}

The bold lines in Table B1 indicate edges that may be subdivided
arbitrarily many times, whereas the dotted lines indicate paths
(possibly degenerate) of arbitrary length.  Thus each
''exceptional graph'' that appears in Table B1 stands not for just
one graph, but rather for a certain countable collection of
graphs.

\begin{Rk}
We note that among partial linear 2-trees, the graphs for which
$M(G) = 2$ are precisely the graphs of two parallel paths.  The
additional graphs for which $M(A) = 2$, the exceptional families
of Table B1, are, of course, partial (not linear) 2-trees with
very special structure.
\end{Rk}

In the previous section, we completely characterized the $C_2$
graphs with $M = 2$. Consider now a general connected partial
2-tree $G$.  If $G$ is a tree, then $M(G)$ is determined by the
path covering number of $G$.  In particular when $G$ is a tree,
$M(G) = 2$ if and only if $G$ is a graph of two parallel paths.
Let us assume then here and throughout this section that $G$
denotes a connected partial 2-tree that is a not a tree. By
sequentially stripping away all degree-one vertices, we arrive at
the maximal induced subgraph of $G$ that is $C_2$, which, by
definition, is the core of $G$. Let us denote this induced
subgraph by ${\rm Core}(G)$.  By inductively applying Lemma
\ref{L4-2}, we see that
\begin{equation*}
M(G) \geq M({\rm Core}(G)),
\end{equation*}
which, combined with Theorem \ref{T4-8}, proves the following lemma:
\begin{La}
Let $G$ be a partial 2-tree.  If $M(G) = 2$, then the core of $G$
is an $LSEAC$ graph and a graph of two parallel paths.
\label{L5-2}
\end{La}
In light of this Lemma,
asking for what connected graphs we have
$M = 2$ is
equivalent to asking under what conditions does
sequentially adding degree one vertices to $LSEAC$ graphs
preserve maximum multiplicity.  Lemma \ref{L4-2} will allow us
to prove several ``forbidden subgraph'' lemmas in which it will
be shown that certain induced subgraphs of a graph $G$ preclude
$M(G) = 2$.  Indeed, if we take an $LSEAC$ graph and sequentially
add degree one vertices until we have created a forbidden subgraph,
then, by Lemma \ref{L4-2} adding subsequent degree one vertices cannot
decrease maximum multiplicity.

\begin{La}
Let $H$ be a $C_2$ graph and let $u$ be an arbitrary vertex of $H$.
Let $T$ be a tree that is a not a path and let $v$ be a degree-one vertex of $T$.
Let $G$ be the graph that is the result of identifying vertices
$u$ and $v$ of $H$ and $T$, respectively.  Then $M(G) > 2$.
\label{L5-3}
\end{La}
\begin{proof} Since $T$ is not a path, it has a vertex of degree at
least three, say, vertex $w$.  We may then find a minimal induced
path connecting $w$ to a pendant vertex of $G$ that also belongs
to the induced subgraph $T$ in $G$.  Let us sequentially
remove the vertices of this path, including vertex $w$.
The resulting graph has at least two connected
components, one of which is
a graph that has a nontrivial core, the other being a tree
(possibly degenerate).
A graph with nontrivial core has $M > 1$, and so the resulting
induced subgraph of two connected components has $M > 2$.
Applying Lemma
\ref{L4-2} inductively shows that this implies that $M(G) > 2$.
\end{proof}

The immediate consequence of this lemma is that if a given graph $G$
cannot be constructed from an $LSEAC$ graph $H$ by sequentially
articulating paths onto vertices of $H$ only, then $M(G) > 2$.
Suppose now we are given a graph $G$ that can be constructed in such a
manner.  Inductively applying Lemma
\ref{L4-2} then shows that in fact we need only consider
the induced subgraph of $G$ obtained by sequentially removing
pendant vertices whose neighbors have degree two, i.e.,
we may assume without loss of generality that all paths
with pendant vertices are of length one.
We shall use this fact implicitly throughout.

\begin{Dn}
A simple partial 2-tree is a partial 2-tree $G$
whose core is a nontrivial $LSEAC$ graph and which may be constructed
from its core through sequential articulation of vertices to
vertices belonging to its core.
\end{Dn}
We consider $LSEAC$ graphs to be simple partial 2-trees.
The remarks immediately preceding the definition show that if
a graph $G$ is not a simple partial 2-tree, then
$M(G) > 2$.
Whenever $v$ is a vertex of $G$ with a neighbor $u$ that is a
pendant vertex in $G$, we shall call $u$ a pendant neighbor of
$v$.

\begin{La}
Let $G$ be a simple partial 2-tree and $u$ a vertex of $G$
with at least three pendant neighbors.  Then $M(G) > 2$.
\label{L5-5}
\end{La}
\begin{proof}
Let $v_1,v_2,v_3$ denote three (distinct) pendant neighbors of $u$.
Then the  induced subgraph $G - \{u,v_3\}$ has three connected
components: isolated vertices $v_1$ and $v_2$, and the induced subgraph
$G - \{u,v_1,v_2\}$.  It follows that
\begin{equation*}
M(G - \{u,v_3\}) = M(v_1) + M(v_2) + M(G - \{u,v_1,v_2\}) \geq 1 + 1 + 1 = 3
\end{equation*}
By Lemma \ref{L4-2}, $M(G) \geq M(G - \{u,v_3\})$, since $v_3$ is
pendant, and therefore $M(G) \geq 3$.
\end{proof}

Next we introduce the definition
of a terminal cycle.  We shall see that these cycles play a rather
more important role than the other cycles of a graph.

\begin{Dn}
A terminal cycle of an $LSEAC$ graph is a cycle that has at most one
neighbor.
\end{Dn}
If $G$ is a graph whose core is an $LSEAC$ graph, then a cycle $Z$ of $G$
is said to be a terminal cycle if and only if $Z$ is a terminal cycle
in the core of $G$.  Note that a single cycle is an $LSEAC$ graph and
is, by our definition, considered to be a terminal cycle.

\begin{La}
Suppose $G$ is a simple partial 2-tree with a cycle $Z$ that is
not terminal.  Suppose further that there exists a vertex $u$ in
$Z$ that is of degree at least two in the core of $G$, does not
belong to a terminal cycle, and has at least one pendant neighbor
$v$ in $G$. Then $M(G)
> 2$. \label{L5-7}
\end{La}
\begin{proof}
First suppose that $u$ is of degree two in the core of $G$.
Note that this requirement is
equivalent to requiring that $u$ belong to cycle $Z$ and no other
cycles.
Observe that there exists at least one vertex of $Z$ that is
a cut vertex in the (nontrivial) $C_2$ graph ${\rm Core}(G - \{u,v\})$.  For example, consider
a vertex $w$ belonging to $Z$ and one of $Z$'s neighbors such that
a minimal path between $w$ and $u$ includes an edge that serves as
an edge of articulation in the core of $G$.
It follows from Lemma
\ref{L4-1} that $M\left({\rm Core}(G - \{u,v\})\right) > 2$.  Applying Lemma \ref{L4-2}
repeatedly,
we get the desired conclusion $M(G) > 2$.

Now suppose that $u$ is of degree at least three in the core
of $G$.  Necessarily, $u$ belongs to at least two cycles, say
$Z_1$ and $Z_2$.  Since $Z_1$ and $Z_2$ have a unique edge
of articulation, say $(u,w)$, it follows that $w$ is a cut
vertex in the graph ${\rm Core}(G - \{u,v\})$.
By Lemma
\ref{L4-1}, $M\left({\rm Core}(G - \{u,v\})\right) > 2$, whence $M(G) > 2$, 
which completes the proof.
\end{proof}

\begin{Dn}
A distinguished vertex of an $LSEAC$ graph $G$ is a vertex that belongs
to every terminal cycle of $G$.
\end{Dn}
If $G$ is a graph whose core is an $LSEAC$ graph, then a vertex $v$
of $G$ is said to be distinguished if and only if $v$ is distinguished
in the core of $G$.  A vertex $u$ of $G$ that is not a distinguished
vertex shall be called nondistinguished.  Note that if $G$ is a single
cycle, then every vertex of $G$ is a distinguished vertex.  If $G$
is an $LSEAC$ graph that is not a cycle, then $G$ has either $0$,
$1$, or $2$ distinguished vertices.  In fact, if $G$ is any graph
whose core is an $LSEAC$ graph with $n$ distinguished vertices,
then $G$ has at least $|n - 3| + 3$ core vertices for $n = 0,1,2$
and, for $n = 3,4,\ldots,$,
exactly $|n - 3| + 3 = n$ core vertices.
In other words, if we suppose that a graph has at least $n$
distinguished vertices, then we have already established at least
a lower bound on the number of vertices that the core of $G$
may have.  On the other hand, an assumption on the number of terminal
cycles $G$
restricts the possibilities for the number of
distinguished vertices.
\begin{Rk}
A vertex of a simple partial 2-tree $G$ is a distinguished vertex
if and only if $G$ has at least one cycle and the vertex belongs
to each cycle of $G$.
\end{Rk}
\begin{proof} Since the core of $G$ is a nontrivial $LSEAC$ graph,
and because $LSEAC$ graphs may be sequentially constructed
by cycle articulation along edges, the statement immediately
follows.
\end{proof}

\begin{La}
If $G$ is an $LSEAC$ graph with at least one distinguished vertex
$v$, then the induced subgraph $G - v$ is a path.
\label{L5-10}
\end{La}
\begin{proof} This is a simple consequence of the fact that a
distinguished vertex belongs to each cycle.  For,
if $G$ is a cycle, then this statement is trivial.
If $G$ has two distinguished vertices, then $G$ must simply
be a pair of cycles articulated along a single edge.  In this case
too it is clear that $G - v$ is a path.  Finally, if $G$ has
one distinguished vertex, then it follows that every edge of
articulation intersects the distinguished vertex, and so only
a path remains in the induced subgraph.
\end{proof}

\begin{La}
Let $G$ be a simple partial 2-tree with at least one
distinguished vertex with precisely
two pendant neighbors. Then $M(G) = 2$ if and only
if $G$ is a graph of two parallel paths.
\label{L5-11}
\end{La}
\begin{proof}
Let $u$ be a distinguished vertex of $G$ with precisely two
pendant neighbors, $v_1, v_2$.
The induced subgraph $G - \{u,v_1\}$ has two connected components:
the isolated vertex $v_2$ and the graph $G - \{u,v_1,v_2\}$.
By Lemma \ref{L4-2}, $M(G) > 2$, unless we have
$M(G - \{u,v_1,v_2\}) = 1$, i.e., unless $G - \{u,v_1,v_2\}$ is
a path.  This occurs if and only if the pendant vertices of $G$
added to the core of $G$  extend the path
${\rm Core}(G) - u$ to a path $P$ in $G$.
In this situation,
the path $P$ and the path connecting $v_1,u,v_2$ constitute
a pair of two parallel paths, showing on account
of Lemma \ref{L3-7} that $M(G) = 2$.
\end{proof}

\begin{La}
Let $G$ be a simple partial 2-tree with a terminal cycle
$Z$ that has two nondistinguished vertices $u_1,u_2$ satisfying the following
conditions: if either vertex is of degree two relative to the core
of $G$, then they are not neighbors; and, $u_1$ and $u_2$ each
have at least one pendant neighbor, say $v_1$ and $v_2$,
respectively.  Then $M(G) > 2$.
\label{L5-12}
\end{La}
\nline

\noindent{\bf Remark.} Note that the conditions imposed on $G$ imply
that the core of $G$ cannot be a cycle,
since a cycle has no nondistinguished vertices.
\nline

\begin{proof}
%
Indeed, by our choice of $u_1$ and $u_2$, there exists a path in
$Z$ connecting $u_1$ and $u_2$ such that, with the possible
exception of $u_1$ and $u_2$, all vertices in the path belong to
cycle $Z$ only.  This path is one of our connected components. Let
us call it $P_1$. The remainder of the graph (excluding
$u_1,u_2,v_1,v_2$) becomes the other connected component, which we
shall call $P_2$. Since $P_2$ contains at least one cycle (e.g.,
the other terminal cycle), it satisfies $M(P_2) > 1$.  Hence
\begin{equation*}
M(G - \{u_1,v_1,u_2,v_2\}) = M(P_1) + M(P_2) > 2.
\end{equation*}
By Lemma \ref{L4-2}, it follows that $M(G) > 2$.
\end{proof}

\begin{La}
Let $G$ be a simple partial 2-tree with a terminal
cycle $Z$ that has a nondistinguished vertex $u$.  Suppose
that $u$ has at least two pendant neighbors, say $v_1$ and
$v_2$.  Then $M(G) > 2$.
\label{L5-13}
\end{La}
\begin{proof} The induced subgraph $G - \{u,v_2\}$ consists of
two connected components: the isolated vertex $v_1$ and the
induced subgraph $G - \{u,v_1,v_2\}$.  Since
$G - \{u,v_1,v_2\}$ has at least one cycle, we see,
by reasoning now entirely analogous to that used in
Lemma \ref{L5-12}, that $M(G) > 2$.
\end{proof}

\begin{La}
Let $G$ be graph whose core is a nontrivial $LSEAC$ graph.
Moreover, suppose that if any distinguished vertex of $G$ has
exactly one pendant neighbor, then there exists another
distinguished vertex of $G$ with at least two pendant neighbors.
Then $M(G) = 2$ if and only if $G$ is a graph of two parallel
paths. \label{L5-14}
\end{La}
\begin{proof} Without loss of
generality, we may assume that $G$ is a simple partial 2-tree
(for if $G$ cannot be constructed from some simple partial 2-tree by
a sequence of subdivisions of edges adjacent to pendant vertices,
then $M(G) > 2$, and so by Lemma \ref{L3-7}, $G$ cannot
be a graph of two parallel paths).
By Lemma \ref{L5-5}, if some distinguished vertex of $G$ has
more than two pendant neighbors, then $M(G) > 2$.  By Lemma
\ref{L3-7}, $G$ is not a graph of two parallel paths.
If some distinguished vertex of $G$ has
precisely two pendant neighbors, then the claim follows by Lemma
\ref{L5-11}.  If no distinguished vertex of $G$ has any
pendant neighbors, then, for $M = 2$, it is necessary that the following be
satisfied: pendant vertices must be adjacent to vertices of
terminal cycles, by Lemma \ref{L5-7}; no vertex of a terminal cycle
may have more than one pendant neighbor, by Lemma \ref{L5-13}; and
no terminal cycle may have more than two nondistinguished vertices
each with at least one pendant neighbor, and moreover,
if a terminal cycle has two such nondistinguished vertices, then those
two vertices
in the terminal cycle are neighbors (Lemma \ref{L5-12}).  When
these conditions are satisfied, however, it follows that there
exist disjoint induced paths $P_1$ and $P_2$ in the core of $G$,
covering all vertices in the core, such that the pendant
vertices of $G$ extend paths $P_1$ and $P_2$.  Hence it is necessary
that $G$ be a graph of two parallel paths.  The sufficiency of
being a graph of two parallel paths has already been shown.
\end{proof}

The following result is a useful criterion for ruling out graphs
for which $M > 2$.
\begin{La}
Let $G$ be a (general) connected graph.  If $G$ has more
than five pendant vertices, then $M(G) > 2$.
\label{L5-15}
\end{La}
\begin{proof}
Without loss of generality, we may assume that $G$ is a simple
partial 2-tree.
If $G$ is a tree with more than five pendant vertices,
then $G$ has path covering number at least three, and so $M(G) \geq 3$.
It therefore
suffices to consider the case where the core of $G$ is a
simple partial 2-tree with at least six pendant vertices.
If no distinguished vertex of $G$ has precisely one pendant neighbor,
 then by Lemma \ref{L5-14}, $M(G) > 2$ since clearly $G$ cannot
be a graph of two parallel paths.  We may therefore restrict ourselves
to the case where some distinguished vertex of $G$, say $u$, has
precisely one pendant neighbor, say $v$.  Since $u$ belongs
to each cycle of $G$, $G - \{u,v\}$ is a tree.
By our assumptions, $G - \{u,v\}$ has at least five
pendant vertices, meaning now that it is a tree with
path covering number at least
three, implying that $M(G - \{u,v\}) > 2$.  Invoking Lemma
\ref{L4-2} one more time, we see that $M(G) > 2$.
\end{proof}

It is a direct consequence of the lemmas we have proven thus far
that in order
for a graph $G$ to be ``exceptional'', i.e., have $M(G) = 2$
and yet not be a graph of two parallel paths, $G$ must also satisfy
the following:
\begin{enumerate}
\item \text{$G$ is a simple partial 2-tree.}
\item \text{$G$ has at least one distinguished vertex.}
\item \text{At least one distinguished vertex has exactly one
pendant neighbor.}
\item \text{No vertex has more than one pendant neighbor.}
\item \text{$G$ has no more than five pendant vertices.}
\end{enumerate}
That each graph belonging to one of the classes of graphs listed in
the table is exceptional may be verified readily: inductively
use the pendant vertex lemma on the path extending from a
distinguished vertex, considering both the case
where the distinguished vertex remains but without the path
attached to it, and the case where the path and the distinguished vertex
are removed.  In each case, the induced subgraph has $M = 2$, and
therefore, by Lemma \ref{L4-2}, so does the original graph.
Clearly, no graph in the table is a graph of two parallel paths.
This is obvious when the number of pendant vertices is five.  When
there are fewer pendant vertices, a small, finite number of subcases
may be considered showing that no graph may be covered by two
disjoint induced paths satisfying the edge crossing condition.

In the subsequent analysis we shall show that these are the
only exceptional graphs.  In proving the results that follow we
shall frequently make use of the fact that the five conditions
listed above are necessary for a graph to be exceptional. In other
words, it suffices to consider only the cases where
the above conditions hold.

\begin{La}
If a graph $G$ is exceptional, then $G$ has more than two pendant
vertices.
\end{La}
\begin{proof}
Let $G$ be a simple partial 2-tree with a distinguished vertex
$u_1$ that has precisely one pendant neighbor, say $v_1$.
Let
$u_2$ be a vertex, different from $u_1$, that belongs to the core of
$G$, and let $v_2$ be a pendant neighbor of $u_2$.  Suppose that
$v_1$ and $v_2$ are the only degree-one vertices of $G$.  The vertex
$u_2$ belongs to a terminal cycle, for otherwise $M(G) > 2$.  In this
case, $u_2$ has a neighbor $w$ (possibly equal to $u_1$)
in the same terminal cycle such that there exists a minimal path
$P_1$ connecting $w$ and $v_1$ such that $G - P_1$ has only
one connected component.  In this case, though, $G - P_1$
is a path, which we denote by $P_2$.  By the construction,
$P_1$ and $P_2$ constitute a pair of parallel paths, showing
that $M(G) = 2$ and $G$ is a graph of two parallel paths.
We have
thus shown that a simple partial 2-tree $G$ with two pendant
vertices has $M(G) = 2$ if and only if $G$ is a graph of two parallel
paths.  This completes the proof.
\end{proof}

\begin{La}
If a graph $G$ is exceptional, then $G$ has no more than three
cycles.
\end{La}
\begin{proof}
Let $G$ be a simple partial 2-tree whose core has at
least four cycles and such that no vertex in the core of $G$
has more than one pendant neighbor.
Since $G$ has at least four cycles, the cycles determine a
unique distinguished vertex,
which we shall call $u$.  Suppose $u$ has a unique pendant
neighbor $v$, and that the graph $G - v$ has either $2$, $3$,
or $4$ pendant vertices and is a graph of two parallel paths.
Note that if $G - v$ is not a graph of two parallel paths,
then, then by Lemma \ref{L5-14}, $M(G) > 2$.
If $G$ is a graph of
two parallel paths, then $G$ is not exceptional, and so we suppose
that $G$ is not a graph of two parallel paths.
In particular, this
implies that in each terminal cycle there is at least one vertex
different from $u$ that has a pendant neighbor.  It also implies
that, if in
terminal cycle $Z$ the neighbor of $u$ that belongs only in cycle
$Z$ supports a pendant vertex, then there are at least
three vertices in $Z$ (including $u$) that have pendant neighbors.
Now consider
the induced subgraph $G - \{u,v\}$.  Since $u$ is distinguished,
the induced subgraph is a tree.  Since $G$ has at least four cycles,
no path shorter than two edges
in length exists between vertices of opposite terminal cycles.
On the other hand, at least
two vertices of each terminal cycle are pendant in the graph
$G - \{u,v\}$ by the conditions we have already established
on pendant neighbors.  Together, these statements
imply that $G - \{u,v\}$ is a tree of path covering number three.
By Lemma \ref{L4-2}, $M(G) > 2$ and so, in particular,
$G$ is not exceptional.
\end{proof}

\begin{La}
A graph $G$ whose core consists of three cycles is exceptional
if and only if $G$ belongs to the collection of graphs
given in Table B1 (see Appendix B).
\end{La}
\begin{proof}
Let $G$ be a simple partial 2-tree whose core consists of
three cycles and let us assume that $G$ has at least three pendant
neighbors. We assume that only vertices of terminal cycles may have
pendant neighbors.
The graph $G$ has a unique distinguished vertex
$u$.  In order for $G$ to be exceptional, $u$ must have precisely
one pendant neighbor, say $v$.  In order for $M(G) = 2$,
$G - \{u,v\}$ must be a tree of path covering number two (if it
were one, then $G$ would be a graph of two parallel paths).  Let
$Z$ be the unique cycle of $G$ that is not a terminal cycle.  If
$Z$ contains more than three vertices, then $M(G - \{u,v\}) > 2$
unless $G$ is a graph of two parallel paths.  Henceforth assume
that $Z$ has only three vertices.  In this case, though, by
considering the graph $G - \{u,v\}$, we see that, in order for
$G$ to be exceptional, each vertex of $Z$ must have precisely
one pendant neighbor.

Suppose we have such a graph $G$, but with the additional
requirement that only the vertices of $Z$ have pendant neighbors.
Then it is easy to check that such a graph $G$ is exceptional and
is in Table B1.  Suppose we wish to add a pendant vertex to $G$ so
that the resulting graph is still exceptional.  Necessarily, we
must add a new vertex to a terminal cycle, say $Z_t$. Considering
$G - \{u,v\}$ shows that it is also necessary that $Z_t$ have only
three vertices.  Considering $G - v$ shows that this is
sufficient.  Such graphs are also in our table. To add a fifth
pendant vertex, necessarily each terminal cycle has only three
vertices, and in that case, there is a unique way pendant vertices
may be arranged.  This case is also exceptional and is covered in
Table B1. Since this analysis exhausts all possibilities, the
proof is complete.
\end{proof}

\begin{La}
A graph $G$ whose core consists of two cycles is exceptional
if and only if $G$ belongs to the collection of graphs given
in Table B1.
\end{La}
\begin{proof}
Let $G$ be a simple partial 2-tree whose core consists of
two cycles and let us assume that $G$ has at least three pendant
neighbors, and that no vertex of $G$ has multiple pendant neighbors.
We first claim that $G$
cannot be exceptional unless each distinguished vertex,
$u_1,u_2$, has precisely
one pendant neighbor.
For suppose that $u_1$ has a unique pendant neighbor $v_1$, but
that $u_2$ does not have a pendant neighbor.  If $G - \{u_1,v_1\}$
is a path, then $G$ is a graph of two parallel paths, and so
suppose that $G - \{u_1,v_1\}$ is a tree of two parallel paths.
Then at least one cycle has at most one nondistinguished vertex
with a pendant neighbor.  However, in this event $G$ must be a graph
of two parallel paths

We therefore may assume that $u_1$ has a pendant neighbor $v_1$ and
$u_2$ has a pendant neighbor $v_2$.
If $M(G) = 2$, then we must have $M(G - \{v_1,v_2\}) = 2$, which
happens only when this graph is a graph of two parallel paths.
In this case we see that if a cycle of $G - \{v_1,v_2\}$ has
two vertices each with a pendant neighbor, then these two
vertices must be adjacent.  If no cycle of this graph has
two vertices, then, by considering $G - \{u_1,v_1\}$, we
see that nondistinguished vertices with pendant neighbors
must be distributed so that one is adjacent to $u_1$ and
the other to $u_2$ (or perhaps both).  In this case, however,
$G$ is a graph of two parallel paths.

Thus we assume that there is a cycle $Z$ of $G$ that has two
nondistinguished adjacent vertices, each with a pendant
neighbor.  By considering the trees $G - \{u_1,v_1\}$
and $G - \{u_2,v_2\}$, we see that $Z$ must have only
four vertices, for otherwiswe $M(G) > 2$.  Let $G$ be a graph
satisfying all of these requirements.  Suppose in addition
that $G$ has only four pendant vertices.  Then $G$ is an
exceptional graph and is in Table B1.  Suppose we wish to
add a fifth pendant vertex to $G$.  This vertex must
necessarily be added to a vertex in $G - Z$.  Considering
$G - \{u_1,v_1\}$ and $G - \{u_2,v_2\}$ shows that this
may only be done if the cycle in $G$ that is not $Z$
is a cycle on three vertices.  The resulting graph is
unique, exceptional, and in our table.  Since we have
exhausted all possibilities, the proof is complete.
\end{proof}

\begin{La}
A graph $G$ whose core is a single cycle is exceptional
if and only if $G$ belongs to the
collection of graphs given in Table B1.
\end{La}
\begin{proof}
Let $G$ be a simple partial 2-tree with between three and five
pendant vertices and such that no vertex has more than one pendant
neighbor.  If the core of $G$ has only three or four vertices,
then clearly $G$ is a graph of two parallel paths.  If the core of
$G$ has five vertices, then $G$ is a graph of two parallel paths
unless $G$ has five pendant vertices.  If $G$ has five vertices in
its core and five pendant vertices arranged as specified above,
then it is exceptional and is in our table.  Now suppose that the
core of $G$ has more than five vertices.  If $G$ has at most two
pendant vertices, then $G$ is a graph of two parallel paths.  Note
that, since $G$ has more than five vertices, there exists a set
(in general, many) of three vertices of the core of $G$, say
$u_1$, $u_2$, $u_3$, such that the induced subgraph $G -
\{u_1,u_2,u_3\}$ has three connected components.  If $u_1$, $u_2$,
and $u_3$ all have a pendant neighbor, then, by Lemma \ref{L4-2}
applied three times, $M(G) > 2$.  If three such vertices, each
with a pendant neighbor, cannot be found, then $G$ has at most
four pendant vertices and is a graph of two parallel paths. This
exhausts all possibilities and so completes the proof.
\end{proof}

\noindent{\bf Remark.} Let $G$ be as in the previous lemma, but with
the additional assumption that the core of $G$ has at least
four vertices.  Suppose $G$ is a graph
of two parallel paths.  If $G$ has three pendant vertices, then
there is at least one pair of pendant vertices that are such that
their neighbors are adjacent.  If $G$ has four pendant
vertices, then we may partition the pendant
vertices into two such pairs.
\nline

\begin{La}
Let $G$ be a simple partial 2-tree with a distinguished
vertex $u$ that has precisely one pendant neighbor, say $v$.
Then $M(G) = 2$ if and only if $G - v$ is a graph of two parallel
paths and $G - \{u,v\}$ is either a graph of two parallel paths
or a path.
\end{La}
\begin{proof}
It is clear that if $G - v$ is a graph of two parallel paths and
$G - \{u,v\}$ is either a graph of two parallel paths or a path,
then $M(G) = 2$.  If $M(G) = 2$ and $G$ is a graph of
two parallel paths, then the implication is trivial.  If $G$
is exceptional, then $G$ is included in our table.  It may be
easily checked that the claim holds for each of these graphs.
\end{proof}
\begin{Rk}
Our main result, Theorem \ref{T5-1}, holds in fact over any
infinite field. This has been explained in detail in the proof of
Lemma \ref{L2-3}, and to a lesser degree in the proof of Lemma
\ref{L4-3}. Many of our lemmas generalize in a straightforward way
from $\mathbb{R}$ to any infinite field (and sometimes to any
field F).
\end {Rk}
\newtheorem{Ob}[Pa]{{\bf Observation}}
\newenvironment{Rem}{\par\noindent{\bf Remark:}\hspace{5mm}}{.\vspace{2mm}\\}
\newcommand {\wT}[1]{\tilde{a}_{#1}}
\newcommand {\hsp}{\hspace{5mm}}
\newpage
\appendix
\renewcommand {\thesection} {Appendix A}
\section {}
\renewcommand {\thesection} {A}
\indent

The purpose of the appendix is to prove Lemma \ref{L2-3}. In fact,
 we prove a more general result, by replacing $\mathbb{R}$ by any
infinite field $F$. The proof given has an algebraic-combinatorial
flavor.

We start with a few preliminaries. Let $S(F,G)$ denote the set of
all $n\times n$ symmetric matrices with entries in $F$, and whose
graph is $G$.  We use $E$ to denote the set of edges of the graph $G$.  Let
\begin{equation*}
\text{msr}(F,G)=\text{min rank}\,A ,
\end{equation*}
where $A$ ranges over all matrices in $S(F,G)$. Let
\begin{equation*}
\text{cork} (F,G)=n-\text{msr}(F,G) .
\end{equation*}
\begin{Rem} We clearly have cork$(\mathbb{R},G)=M(G)$ \end{Rem}

Let $H$ be any subgraph of $G$ which is an $hK_4$. We call the
original vertices of $K_4$, used to obtain $H$ from $K_4$ by a
sequence of edge subdivisions, {\em initial vertices}. All other
vertices of $H$ are called {\em intermediate vertices}.\nline

The following well known result is used in the Appendix:
\begin{Ob}
Let $F$ be an infinite field and let $f\in F[t_1,t_2,\cdots,t_n]$.
Then there exist $a_1,a_2,\cdots,a_n\in F^*$ such that
$f(a_1,a_2,\cdots,a_n)\in F^* .$
\end{Ob}

The following discussion plays an important role in the proof of
the main result of the Appendix.

Let $l$ be an integer such that $3\leq l\leq n-1$. Let $G'$ be the
induced subgraph on $\{l+1,l+2,\cdots,n\}$. Let $B$ be an
$n-l\times n-l$ symmetric matrix defined as follows: its rows and
columns are labeled from $l+1$ to $n$ ; its main diagonal entries
are all zero. For $i,j \in \{l+1,l+2,\cdots,n\}, i<j$, we let
$b_{ji}=b_{ij}=0$ if $ij \notin E$ and we let $b_{ij}$ be an
indeterminate if $ij \in E$ (and let $b_{ji}=b_{ij}$). So, for
example, if $l+1,l+2\in E$ then the entry in the first row and
second column of B is $b_{l+1,l+2}$. We let $\underline{b}$ denote
the set of all indeterminates in $B_{22}$.

\noindent Let
\begin{equation*} A_{22}=I-\mu B_{22} ,\end{equation*}
where $\mu$ is an indeterminate. Let $A_{12}=(a_{ij}) ,\hsp
i=1,2,\cdots,l \hsp  j=l+1,l+2,\cdots,n$ be such that the $ij$-th
entry is 0 if $ij \notin E$ and an indeterminate $a_{ij}$ if $ij
\in E$.

Let $R$ be the polynomial ring consisting of all polynomials in
the indeterminates that appear in $A_{12}$ and $A_{22}$, and with
coefficients in $F$. Let $K$ be the quotient field of $R$. Let
\begin{equation*} d=d(\mu,\underline{b})=\text{det} A_{22} ,\end{equation*}
\begin{equation*} W=W(\mu,\underline{b})=\text{adj} A_{22} ,\end{equation*}
and
\begin{equation*} Z=Z(\mu,\underline{b})=A_{22}^{-1} .\end{equation*}
Note that $d \in K^*$, because the constant term in its expansion
is 1. Hence $A_{22}$ is an invertible element of $K^{n-l,n-l}$ ,
so $Z$ exists and we have
\begin{equation*} W=dZ\in K^{n-l,n-l} .\end{equation*}

\begin{Rem}
As for $B_{22}$, the rows and columns of $A_{22}, W$ and $Z$ are
labeled $l+1$ to $n$
\end{Rem}
\noindent Define
\begin{equation*} C=A_{12}A_{22}^{-1}A_{12}^t=A_{12}ZA_{12}^t \in K^{l,l} ,\end{equation*}
and let q be a positive integer, to be determined later. We have
\begin{equation*} (I-\mu B_{22})(I+\mu B_{22}+\mu ^2B_{22}^2+\cdots +\mu ^qB_{22}^q)=I-\mu^{q+1}B_{22}^{q+1} , \end{equation*}
hence
\begin{equation*} I+\mu B_{22}+\mu ^2B_{22}^2+\cdots +\mu ^qB_{22}^q=Z-\mu ^{q+1}ZB_{22}^{q+1} . \end{equation*}
Since $Z=d^{-1}W$, we get
\begin{equation} Z=d^{-1}[d(I+\mu B_{22}+\mu ^2B_{22}^2+\cdots +\mu ^qB_{22}^q)+\mu ^{q+1}WB_{22}^{q+1}] , \label{A-1} \end{equation}
and
\begin{equation} W=d(I+\mu B_{22}+\mu ^2B_{22}^2+\cdots +\mu ^qB_{22}^q)+\mu ^{q+1}WB_{22}^{q+1} . \label{A-2} \end{equation}
\\
\begin{La}\label{LA-2}
Let $i,j\in \{1,2,...,l\}$ with $i\neq j$. Suppose that there
exists a path in $G \hsp ii_1,i_1i_2,...,i_ri_{r+1},i_{r+1}j$ such
that $i_1,i_2,...,i_r,i_{r+1}\notin \{1,2,...,l\}$. Then
$c_{ij}\in K^*$.
\end{La}
\begin {proof}
We show first that $z_{i_1,i_{r+1}}\neq 0$. It suffices to show
$w_{i_1,i_{r+1}}\neq 0$. Indeed, in computing $w_{i_1,i_{r+1}}$
the expression $\mu^r\displaystyle\prod_{k=1}^r b_{i_k,i_{k+1}}$
appears as one of the summands of $(\mu^rB_{22}^r)_{i_1,i_{r+1}}$.
This term cannot appear in any matrix that comes after $\mu^r
B_{22}^r$ in (\ref{A-2}) (we assume $q$ large enough, say
$q>n-l$). It cannot come from preceding matrices either. Indeed,
if it would, some contribution $\neq 1$ must come from $d$. But in
every summand $\neq 1$ in the expression of d the indeterminates
from $\underline{b}$ that appear are $a$ disjoint union of several
cycles, at least one not the identity, and this is impossible
here. We notice now that one of the terms in $dc_{ij}$ is
$a_{i,i_1}a_{i_{r+1},j}\mu^r\displaystyle\prod_{k=1}^r
b_{i_k,i_{k+1}}$ , and it cannot be canceled. Hence $dc_{ij} \in
K^*$, implying that $c_{ij} \in K^*$.
\end {proof}

The main result of the appendix is:
\begin {La} Let $F$ be an infinite field and let $G$ be a graph on
$n \geq 4$ vertices, which is not a partial 2-tree. Then
\begin{equation*} \text{cork}(F,G)\geq3 ,\end{equation*}
\end {La}
\begin {proof}
By Lemma \ref {L2-2} G contains (as a subgraph) an $hK_4$ , which
is denoted by $H$, and its initial vertices are denoted by
$1,2,3,n$. We distinguish four cases :

\Cs There are three initial vertices, say $1,2,3$, such that the
three paths in $H$ contain
no intermediate vertices.\\
\centerline{\includegraphics[width=70mm]{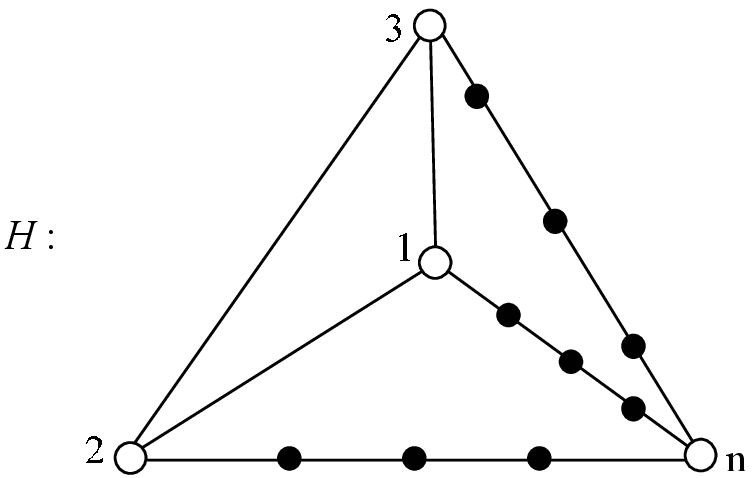}}
\begin{Rem}
Here, and in subsequent figures,\space\space{\Large$\bullet\;$}
denotes a vertex that might or might not be present in $H$
\end{Rem}
\noindent We apply Lemma \ref{LA-2} with $l=3$. Since F is
infinite we can clearly assign nonzero values to all
indeterminates in $K$ so that all off-diagonal entries of
$C=A_{12}A_{22}^{-1}A_{12}^t$ are nonzero. Then the matrix
$\begin{bmatrix}
C & A_{12}\\
A_{12}^t & A_{22}\\
\end{bmatrix}$
is in $S(F,G)$ and has rank $n-3$.%
\Cs There are three initial vertices, say $1,2,3$, such that
exactly one of the three corresponding paths in $H$ contains an
intermediate vertex.\\
\centerline{\includegraphics[width=70mm]{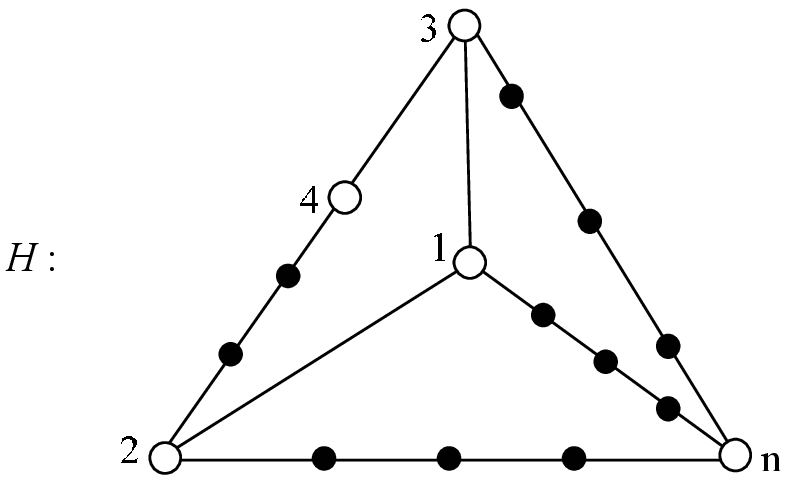}}\\
We assume that 3 and 4 are adjacent, but 2 and 4 \underline{don't
have} to be adjacent.\\
If $24\in E$ we may assume it is also an edge of $H$, for if it
isn't, we replace the path from 2 to 4 in $H$ by the edge $24 $.
Given $i,j \in \{1,2,3,4\}$ such that $i\neq j$, we say there is
an {\em external path} from $i$ to $j$ in $G$ if there is a
(simple) path from $i$ to $j$, which is not an edge, and such that
no intermediate vertices of this path belong to $\{1,2,3,4\}$. In
our case there are external paths from 1 to 2; from 1 to 3 ; from
2 to 3 ; from 2 to 4 if they are not adjacent.\\
We want to use Lemma \ref{LA-2} with $l=4$. Let
\begin{equation*}
A_{11}=A[1,2,3,4]= \begin{bmatrix}
? & ? & ? & a_{14}\\
? & ? & 0 & a_{24}\\
? & 0 & ? & ?\\
a_{14} & a_{24} & ? & ?\\
\end{bmatrix},
\end {equation*}
where we use the following notation: '?' means an indeterminate
(which is NOT one of those that appear in the definition of K);
$a_{14}=0$ if $ 14 \notin E$ , and is chosen in $F^*$ if $14 \in
E$ ; put $a_{24}=0$ if $24\notin E$ and a '?' if $24 \in E$.
\begin{Rem}
We have $a_{23}=0$ because we may assume $23 \notin E $ (or else
apply case 1)
\end{Rem}\noindent
We assign nonzero values to all indeterminates in $K$ so that
$A_{22}$ is nonsingular and every $c_{ij}\in K^*$ becomes an
element of $F^*$. Now, the Schur complement with respect to
$A_{22}$ is
\begin{equation*}
A_{11}-C=A_{11}-A_{12}A_{22}^{-1}A_{12}^t\;.
\end{equation*}
Also, if $a_{24}=?$, we choose it so that the element in the $2,4$
position of $A_{11}-C$ is nonzero. So $A_{11}-C$ has the form:

\begin{equation*}
\begin{bmatrix}
? & ? & ? & \wT{{14}}\\
? & ? & \wT{{23}} & \wT{{24}}\\
? & \wT{{23}} & ? & ?\\
\wT{{14}} & \wT{{24}} & ? & ?\\
\end{bmatrix},
\end{equation*}
where we know $\wT{{23}} , \wT{{24}}\in F^*$. Now let
$x^t=(\dfrac{\wT{{14}}}{\wT{{24}}},1,\wT{{23}},\wT{{24}})$. Then
\begin{equation*}
xx^t=
\begin{bmatrix}
\vspace{2mm} \frac{\wT{{14}}^2}{\wT{{24}}^2} &
\frac{\wT{{14}}}{\wT{{24}}} & \frac{\wT{{14}}\wT{{23}}}{\wT{{24}}}
& \wT{{14}}\\ \vspace{2mm} \frac{\wT{{14}}}{\wT{{24}}} & 1 &
\wT{{23}} & \wT{{24}}\\ \vspace{2mm}
\frac{\wT{{14}}\wT{{23}}}{\wT{{24}}} & \wT{{23}} & \wT{{23}}^2 &
\wT{{23}}\wT{{24}}\\
\wT{{14}} & \wT{{24}} & \wT{{23}}\wT{{24}} & \wT{{24}}^2\\
\end{bmatrix},
\end{equation*}
and it has the desired form. So we can choose $A\in S(F,G)$ with
$\text{rank} A =n-4+1=n-3$.%
\Cs There are three initial vertices, say $1,2,3$, such that
exactly two of the three corresponding paths in $H$ contain an
intermediate vertex.\\
\centerline{\includegraphics[width=70mm]{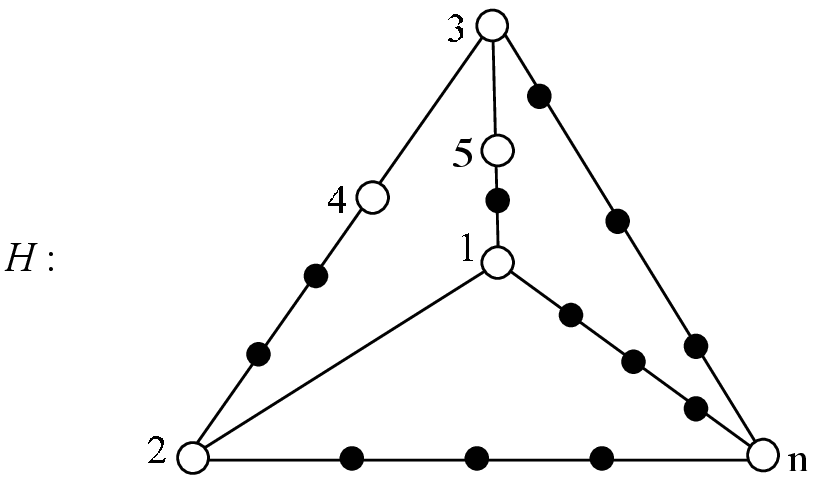}}\\
Here we assume that $34\in E, 35 \in E$ but $24$ and $15$ don't
have to be edges of G. As in Case 2, if $24\in E$ we may assume it
is an edge of $H$. A similar statement for $15$. We use the notion
external path as in the previous case, so, for example, there is
an external path from 1
to 2. If $24 \notin E$, there is an external path from 2 to 4.\\
We want to use Lemma \ref{LA-2} with $l=5$. Let
\begin{equation*}
A_{11}=A[1,2,3,4,5]= \begin{bmatrix}
? & ? & 0 & a_{14} & a_{15}\\
? & ? & 0 & a_{24} & a_{25}\\
0 & 0 & ? & ? & ?\\
a_{14} & a_{24} & ? & ? & a_{45}\\
a_{15} & a_{25} & ? & a_{45} & ?\\
\end{bmatrix}\, ,
\end {equation*}
where '?' is used as before; $a_{14}=0$ if $14 \notin E$ and
$a_{14}\in F^*$ if $14 \in E$; put $a_{15}=0$ if $15\notin E$ and
a '?' if $15 \in E$; put $a_{24}=0$ if $24\notin E$ and a '?' if
$24 \in E$ ; $a_{25}=0$ if $25 \notin E$ and $a_{25}\in F^*$ if
$25 \in E$; $a_{45}=0$ if $45\notin E$ and $a_{45}\in F^*$ if $45
\in E$. Note that we may assume $13 \notin E$ and $23 \notin E$ ,
or else we can apply a previous case. Assign nonzero values to all
indeterminates in K so that $A_{22}$ is invertible and every
$c_{ij}\in K^*$ becomes an element of $F^*$. Now, the Schur
complement with respect to $A_{22}$ is
\begin{equation*}
A_{11}-C=A_{11}-A_{12}A_{22}^{-1}A_{12}^t\;.
\end{equation*}
Also, if $a_{15}=?$ (resp. $a_{24}=?$), we assign it a value so
that $1,5$ entry (resp. $2,4$ entry) of the Schur complement is
nonzero. Hence, $A_{11}-C$ has the form:

\begin{equation*}
\begin{bmatrix}
? & ? & \wT{{13}} & \wT{{14}} & \wT{{15}}\\
? & ? & \wT{{23}} & \wT{{24}} & \wT{{25}}\\
\wT{{13}} & \wT{{23}} & ? & ? & ?\\
\wT{{14}} & \wT{{24}} & ? & ? & \wT{{45}}\\
\wT{{15}} & \wT{{25}} & ? & \wT{{45}} & ?\\
\end{bmatrix}\, ,
\end {equation*}
where $\wT{{13}},\wT{{23}},\wT{{15}},\wT{{24}}\in F^*$ . One needs
to assign values to the indeterminates so that
$A_{11}-C$ has rank two. \\
The $2 \times 2$ principal submatrix in the top left corner is
arbitrary, so we can assume it is invertible and its inverse is
also an arbitrary invertible matrix. Denote it by$\begin{bmatrix}
u & v\\
v & w\\
\end{bmatrix}$. Let $B$ denote its Schur complement. We need
$B=0$. Because of the form of $(A_{11}-C)[3,4,5]$ , we get
effectively only one equation, namely:
\begin{equation*}
0=b_{23}=\wT{{45}}-[\wT{{14}}\;\wT{{24}}]%
\begin{bmatrix}
u & v\\
v & w\\
\end{bmatrix}
\begin{bmatrix}
\wT{{15}}\\
\wT{{25}}\\
\end{bmatrix}
=\wT{{45}}-\wT{{14}}\wT{{15}}u-(\wT{{14}}\wT{{25}}+\wT{{15}}\wT{{24}})v-\wT{{24}}\wT{{25}}w\;
.
\end{equation*}

Suppose first that $\wT{{14}}\wT{{25}}+\wT{{15}}\wT{{24}}\neq 0$.
 If $\wT{{45}} \neq 0$ we have a solution with $u=w=0$ and $v \neq
0$ . So suppose $\wT{{45}}=0$. If $\wT{{14}}=\wT{{25}}=0$ we have
a solution with $v=0$, and $u \neq0 , w \neq 0$. If exactly one of
$\wT{{14}},\wT{{25}}$ is zero , we have a solution with $v \neq 0$
and exactly one of $u,v$ is zero. So we may assume now
$\wT{{14}}\wT{{25}}+\wT{{15}}\wT{{24}}=0$, implying
$\wT{{14}}\wT{{25}}=-\wT{{15}}\wT{{24}}\neq 0\, .$
It is clear that we can find $u,v,w \in F$ such that
$\begin{bmatrix}
u & v\\
v & w\\
\end{bmatrix}$ is invertible. We conclude that there exists $A \in S(F,G)$ with $\text{rank} A=n-5+2=n-3 .$
\Cs The remaining case to consider is of the form:\\
\centerline{\includegraphics[width=70mm]{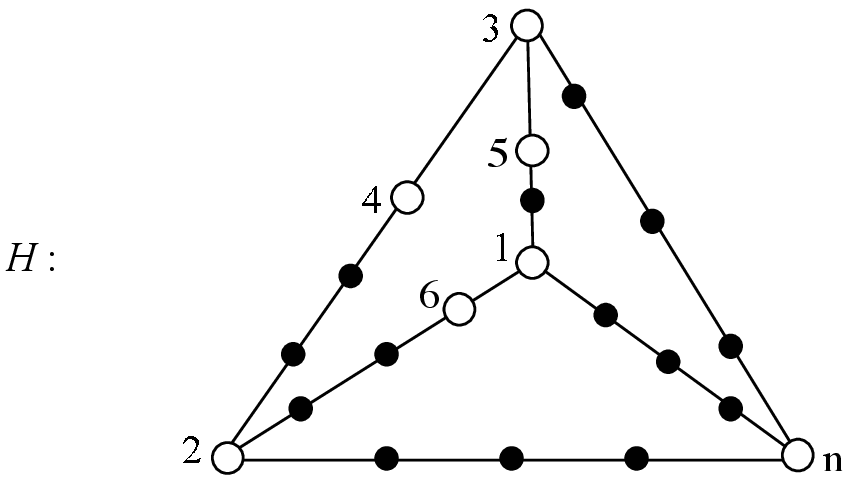}}\\
where $34 \in E,35 \in E ,16 \in E$ , while $24,26$ and $15$ don't
have to be edges. As in previous cases, if any of $24,26$ and $15$
are in $E$ we may assume that they are in $H$. The discussion is
similiar to the previous cases and leads to
\begin{equation*}
A_{11}-C=
\begin{bmatrix}
? & \wT{{12}} & \wT{{13}} & \wT{{14}} & \wT{{15}} & ?\\
\wT{{12}} & ? & \wT{{23}} & \wT{{24}} & \wT{{25}} & \wT{{26}}\\
\wT{{13}} & \wT{{23}} & ? & ? & ? & \wT{{36}}\\
\wT{{14}} & \wT{{24}} & ? & ? & \wT{{45}} & \wT{{46}}\\
\wT{{15}} & \wT{{25}} & ? & \wT{{45}} & ? & \wT{{56}}\\
? & \wT{{26}} & \wT{{36}} & \wT{{46}} & \wT{{56}} & ?\\
\end{bmatrix}\, ,
\end{equation*}
where
$\wT{{12}},\wT{{13}},\wT{{15}},\wT{{23}},\wT{{24}},\wT{{26}}\in
F^*$. One needs to assign values to the indeterminates so that
$A_{11}-C$ has rank 3. We compute the Schur complement $B$ with
respect to the $3 \times 3$ principal submatrix of $A_{11}-C$
based on raws $1,2,6$. We write this principal submatrix as
\begin{equation*}G=
\begin{bmatrix}
x & \wT{{12}} & u\\
\wT{{12}} & y & \wT{{26}}\\
u & \wT{{26}} & z\\
\end{bmatrix}\, ,
\end{equation*}
where $x,y,u,z$ are indeterminates. So
\begin{equation*}
B=
\begin{bmatrix}
? & ? & ?\\
? & ? & \wT{{45}}\\
? & \wT{{45}} & ?\\
\end{bmatrix}-
\begin{bmatrix}
\wT{{13}} & \wT{{23}} & \wT{{36}}\\
\wT{{14}} & \wT{{24}} & \wT{{46}}\\
\wT{{15}} & \wT{{25}} & \wT{{56}}\\
\end{bmatrix}G^{-1}
\begin{bmatrix}
\wT{{13}} & \wT{{14}} & \wT{{15}}\\
\wT{{23}} & \wT{{24}} & \wT{{25}}\\
\wT{{36}} & \wT{{46}} & \wT{{56}}\\
\end{bmatrix}=0 .
\end{equation*}
Because of the special form of the first (matrix) summand, this
leads to just one scalar equation, namely,
\begin{equation*}
\wT{{45}}-
\begin{bmatrix}
\wT{{14}} &  \wT{{24}} & \wT{{46}}\\
\end{bmatrix}
G^{-1}
\begin{bmatrix}
\wT{{15}}\\ \wT{{25}} \\ \wT{{56}}\\
\end{bmatrix}=0 .
\end{equation*}
This is equivalent to : $G$ is invertible and
\begin{equation}
\wT{{45}}\text{det}G-
\begin{bmatrix}
\wT{{14}} &  \wT{{24}} & \wT{{46}}\\
\end{bmatrix}
\text{adj}G
\begin{bmatrix}
\wT{{15}}\\ \wT{{25}} \\ \wT{{56}}\\
\end{bmatrix}=0 .
\label{A-3}
\end{equation}
We have
\begin{equation*}
\text{adj}G=
\begin{bmatrix}
yz-\wT{{26}}^2 & \wT{{26}}u-\wT{{12}}z & \wT{{12}}\wT{{26}}-uy\\
\wT{{26}}u-\wT{{12}}z & xz-u^2 &
\wT{{12}}u-\wT{{26}}x\\
\wT{{12}}\wT{{26}}-uy & \wT{{12}}u-\wT{{26}}x &
xy-\wT{{12}}^2\\
\end{bmatrix}\, .
\end{equation*}
We consider the left hand side of (\ref{A-3}) as a linear function
of $x$, with coefficients in $F[y,z,u]$ . Denote by $\psi (y,z) $
and $-\varphi(y,z,u)$, respectively, the coefficient of $x$ and
the free coefficient. Then \begin{equation*} \psi (y,z)
=\wT{{45}}(yz-\wT{{26}}^2)-\wT{{24}}\wT{{25}}z-\wT{{46}}\wT{{56}}y+\wT{{26}}(\wT{{24}}\wT{{56}}+\wT{{46}}\wT{{25}})
\, .
\end{equation*}
Suppose first that $\psi (y,z)=0 $. This implies
$\wT{{45}}=\wT{{25}}=\wT{{56}}=0$, and hence the constant
coefficient in the right hand side of (\ref{A-3}) is
\begin{equation*}
-\wT{14}\wT{15}(yz-\wT{26}^2)-\wT{24}\wT{15}(\wT{26}u-\wT{12}z)-\wT{46}\wT{15}(\wT{12}\wT{26}-uy)\,
.
\end{equation*}
The coefficient of u in this expression is
$-\wT{15}\wT{24}\wT{26}+\wT{15}\wT{46}y$. We choose $y$ and $z$ in
$F$ so that this coefficient of $u$ is nonzero and $yz-\wT{26}^2
\neq 0$ . We can determine $u$ so that (\ref{A-3}) is satisfied,
and then we determine $x$ so that $\text{det}G\neq 0$.

We assume now that $\psi (y,z) \neq 0$. We let $x=\frac{\varphi
(y,z,u)}{\psi (y,z)}$ , so any choice of $y,z,u$ in $F$ such that
$\psi (y,z) \neq 0$ will yield a solution of (\ref{A-3}). We have
to make a choice such that $\text{det}G \neq 0$. We have
\begin{equation*}
\text{det}G=\frac{(yz-\wT{26}^2)\varphi (y,z,u)}{\psi
(y,z)}+2\wT{12}\wT{26}u-u^2y-\wT{12}^2z .
\end{equation*}
Let
\begin{equation}
\label{A-4} p(y,z,u)=(yz-\wT{26}^2)\varphi
(y,z,u)+(2\wT{12}\wT{26}u-u^2y-\wT{12}^2z)\psi (y,z) \, .
\end{equation}
If $\wT{25} \neq 0$ the coefficient of $uz$ in $p(y,z,u)$ is
$-2\wT{12}\wT{26}\wT{24}\wT{25}\neq 0$, so $p(y,z,u) \neq 0$.
Hence we may assume $\wT{25}=0$. If $\wT{45} \neq 0$ then the
coefficient of $u^2y^2z$ in $p(y,z,u)$ is $\wT{45}$, so
$p(y,z,u)\neq 0$. Hence we may assume $\wT{45}= 0$. Since $\psi
(y,z)\neq 0$ we must have $\wT{56}\neq 0$, implying that the
coefficient of $u^2y$ in $p(y,z,u)$ is $-\wT{24}\wT{26}\wT{56}\neq
0$, so $p(y,z,u)\neq 0$. It follows that $x,y,z,u$ can be chosen
so that $G$ is invertible and (\ref{A-3}) holds, so $B=0$
\end {proof}

\newpage
\appendix
\renewcommand {\thesection} {Appendix B}
\section {}
\renewcommand {\thesection} {A}
\indent

\includegraphics[width=.9\textwidth]{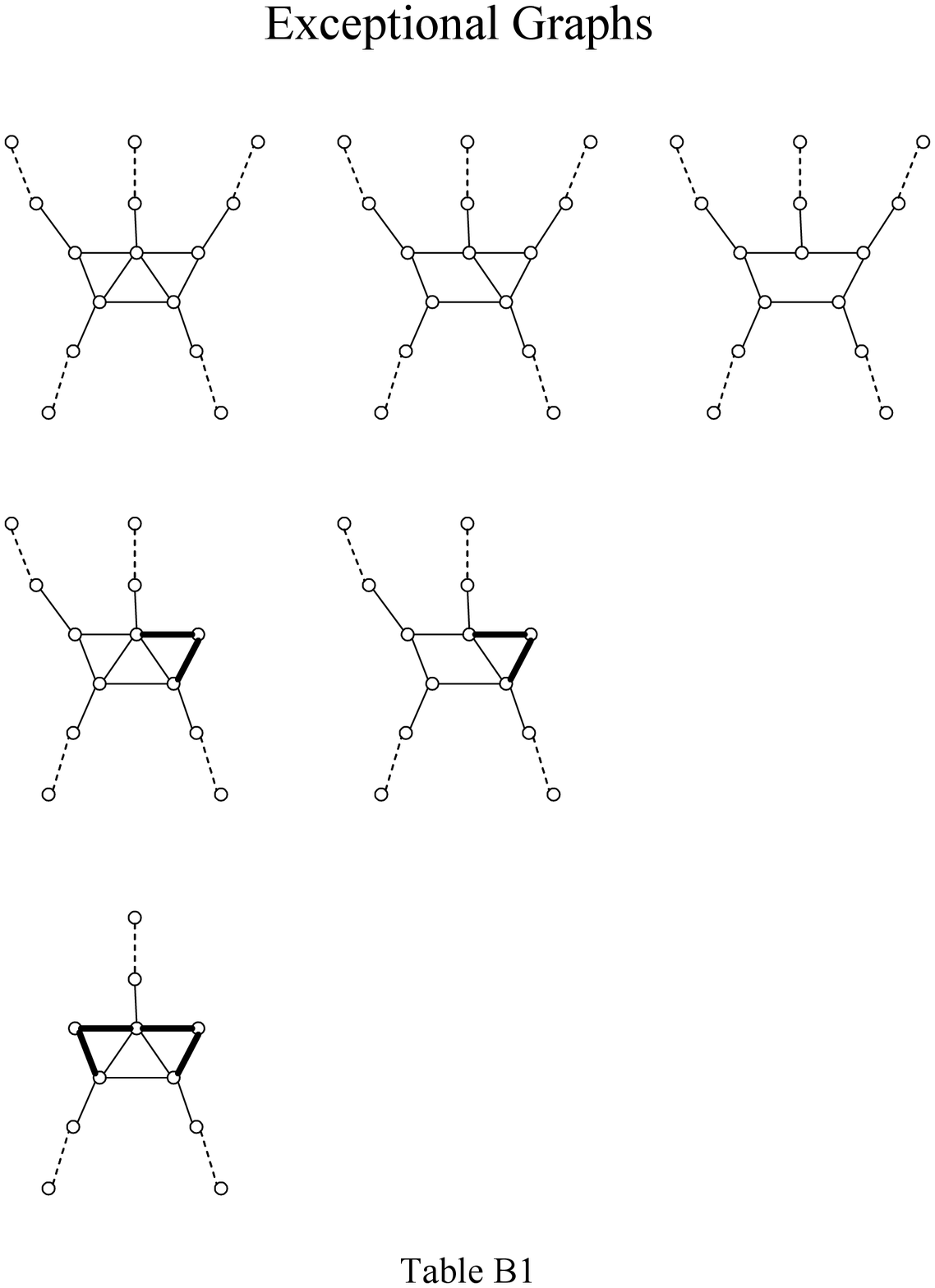}

\end{document}